\documentclass[10pt,a4paper]{article}
\makeatletter
\usepackage[dvipsnames]{xcolor}
\usepackage{lipsum} 
\usepackage{tipa}
\usepackage{tikz}
\usepackage{tikz-cd}
\usetikzlibrary{matrix,arrows,decorations.pathmorphing}
\newcounter{magicrownumbers}

\usepackage[left=3cm,top=3cm,bottom=3cm,right=3cm,head=1cm,foot=1cm]{geometry}
\usepackage[pdftex,bookmarks=true,bookmarksnumbered=true]{hyperref}
\hypersetup{linkbordercolor=ForestGreen!25, citebordercolor=Goldenrod!60}
\DeclareMathAlphabet{\mathpzc}{OT1}{pzc}{m}{it}
\usepackage{graphicx}
\usepackage{stmaryrd}
\usepackage{twoopt}
\usepackage{amssymb}
\usepackage{amsmath}
\usepackage{amsthm}
\usepackage{mathabx} 
\usepackage{amsfonts}
\usepackage{amscd}
\usepackage{url}
\usepackage[all]{xy}

\newcommand{\hil}{\mathcal{H}}

\newcommand{\5}{\hspace{0,5cm}}
\newcommand{\3}{\vspace{0,3cm}}

\newcommand{\ra}{\rangle}

\title{Dirac Operators on Quantum Weighted Projective Spaces}
\date{}
\author{Antti J. Harju\footnote{Copernicus Center for Interdisciplinary Studies, Krakow, Poland} \footnote{harjuaj@gmail.com}}
\begin{document}

\maketitle

\begin{abstract}
The quantum weighted projective algebras $\mathbb{C}[\mathbb{WP}_{k,l,q}]$ are coinvariant subalgebras of the quantum group algebra $\mathbb{C}[SU_{2,q}]$. For each pair of indices $k,l$, two $2$-summable spectral triples will be constructed. The first one is an odd spectral triple based on coinvariant spinors on $\mathbb{C}[SU_{2,q}]$. The second one is an even spectral triple. \3


\end{abstract}

\section*{Introduction}

The theory of the quantum weighted projective algebras incorporates ideas from representation theory, noncommutative geometry and the geometry of singular spaces. The family of algebras $\mathbb{C}[\mathbb{WP}_{k,l,q}]$ is parametrized by two positive coprime integers, $k$ and $l$, and each one is a coinvariant subalgebra in the quantum group $\mathbb{C}[SU_{2,q}]$ under a coaction of the Hopf algebra $\mathbb{C}[u, u^{-1}]$, \cite{BF12}. The Hopf algebra $\mathbb{C}[u, u^{-1}]$ can be identified with the coordinate algebra of a unit circle. Especially interesting cases occur when $k = 1$ and $l$ is a positive integer greater than one when these algebras are $q$-deformations of the coordinate algebras of  singular spaces with a teardrop shape. In the classical case these have been studied extensively in groupoid theory. In this context, a singular space can be naturally associated with a Morita equivalence class of a groupoid \cite{MM03}, \cite{Moe02}. Morita equivalence preserves the shape of the singular space and the types of the singularities but there are several ways of how one can construct the geometric realization of the groupoid. For example, in the case of a teardrop, a representative of the Morita class is usually constructed as follows: there is a collection of open balls and one of them is subject to an action of a finite rotation group. The orbit space of this action is a cone. The balls can be glued together using groupoid arrows so that the orbit space of the groupoid has a shape of a teardrop. The singularity is determined by the choice of the finite rotation group: the parameter $l$ above corresponds to the isotropy at the singularity under the action of $\mathbb{Z}_l$. Another realization for the Morita class is given by an action groupoid determined by an action of the Lie group $\mathbb{T}$ on the sphere $S^3$. The algebraic deformation $\mathbb{C}[\mathbb{WP}_{k,l,q}]$ is based on this model. One can also realize the teardrop spaces as 2-dimensional orbifolds, \cite{Sat56}. The teardrops are examples of 2-dimensional orbifolds which are not global quotients under a finite group action. 

The goal of this work is to put the algebras $\mathbb{C}[\mathbb{WP}_{k,l,q}]$ into the framework of Connes' noncommutative geometry. The coaction of a Hopf algebra on a noncommutative algebra is a noncommutative geometric analogue of a manifold with a group action, or an action groupoid, while the coinvariant subalgebra models the quotient space under the action. In Lie groupoid theory, the local group actions give  rise to local diffeomorphisms. These local diffeomorphisms are applied to define a groupoid action on the tangent bundle. Then one can try to lift the action to the spinor bundle. This is not generally possible, and the obstruction is measured by the Stiefel-Whitney classes in the Lie groupoid cohomology. Whenever the lift exists one can proceed to define an invariant Dirac operator acting on the spinors. This data can be used to model the quotient space as a spectral triple consisting of the invariant function algebra, the Hilbert space of invariant spinors and the invariant Dirac operator, \cite{Har13a} \cite{RV08}. 

Two approaches to define a Hilbert space of spinors and a Dirac operator are taken here, both are based on the representation theory of the Drinfeld-Jimbo algebra $U_q(\mathfrak{su}_2)$, \cite{Dri86, Jim85}. In terms of representation theory, the algebra $U_q(\mathfrak{su}_2)$ has the same complex irreducible representations as the Lie algebra $\mathfrak{su}_2$, but the tensor category of representations has a nontrivial braiding, see e.g. \cite{NT11}. This results in a violation of the permutation symmetry in the representation category which is seen as a noncocommutativity of the coproduct in $U_q(\mathfrak{su}_2)$. In the spirit of Peter-Weyl theorem, one can define the coordinate algebra  $\mathbb{C}[SU_{2,q}]$ as a Hopf dual of $U_q(\mathfrak{su}_2)$ which is linearly spanned by the matrix elements of the irreducible representations, \cite{Wor87b}. The tensor structure in the category of representations determines a product in $\mathbb{C}[SU_{2,q}]$. The product is noncommutative since the coproduct of  $U_q(\mathfrak{su}_2)$ is noncocommutative. The algebra $U_q(\mathfrak{su}_2)$ is given a right and a left representation on the space of matrix elements $\mathbb{C}[SU_{2,q}]$. These representations correspond to the actions of the invariant vector fields on the space of functions. We make the algebra $U_q(\mathfrak{su}_2)$ a $\mathbb{C}[u,u^{-1}]$ comodule algebra by requiring a compatibility for the coaction with the right $U_q(\mathfrak{su}_2)$ representation. If $M_{\frac{1}{2}}$ denotes a $U_q(\mathfrak{su}_2)$ module of the highest weight $\frac{1}{2}$, then the Hilbert space of spinors is naturally defined to be a completion of $\mathbb{C}[SU_{2,q}] \otimes M_{\frac{1}{2}}$. The right $U_q(\mathfrak{su}_2)$ representation extends to this space. The Hilbert space can now be equipped with a comodule structure which is compatible with the right representation. This is a quantum group analogue of a groupoid action on the spinors. Then we can define the coinvariant subspace. This subspace is invariant under the action of the Dirac operator on $SU_{2,q}$. This gives rise to an odd coinvariant spectral triple. The coinvariant spectral triple does not have a chiral grading. This is somewhat expected since it is based on a construction of an odd spectral triple in dimension 3. In the case of a teardrop, $k =1$ and $l \in \mathbb{N}$, another odd spectral triple based on the coinvariance has been developed in \cite{SV13}.

In the second construction we take the coaction on the spinor module $M_{\frac{1}{2}}$ to be trivial. In this case the Hilbert space of spinors is a sum of two copies of the Hilbert space completion of the coinvariant subspace $\mathbb{C}[\mathbb{WP}_{k,l,q}]$. The Dirac spectrum is fixed by requiring that for $k = l = 1$ it gives the classical Dirac spectrum on the manifold $S^2$. The chiral grading can be defined and we obtain an even spectral triple.

The algebras $\mathbb{C}[\mathbb{WP}_{k,l,q}]$ have $C^*$-algebra completions and there are isomorphisms $C(\mathbb{WP}_{k,l,q}) \simeq C(\mathbb{WP}_{k',l,q})$ for all $k$ and $k'$. The geometric models based on the algebras $\mathbb{C}[\mathbb{WP}_{k,l,q}]$ and $\mathbb{C}[\mathbb{WP}_{k',l,q}]$ are never unitary equivalent if $k \neq k'$. Especially the multiplicities of the Dirac eigenvalues depend on the parameter $k+l$ and therefore the Dirac operators associated to different $k$ but equal $l$ cannot be unitary conjugates of each other. Therefore these algebras are equivalent in the topological sense while they have different geometric structures. 

In the course of analysis of the spectral triples we study the decomposition of the space $\mathbb{C}[SU_{2,q}]$ into homogeneous components under the coaction of $\mathbb{C}[u,u^{-1}]$. In the case of quantum teardrops, we can give a full analytic realization for the K-theory groups of $C(\mathbb{WP}_{1,l,q})$ in terms of the completions of these  homogeneous components. The analysis relies and generalizes the recent results based on the study of certain quantum principal $\mathbb{T}$-bundles over the quantum teardrops, \cite{She14}. 

The main reference is \cite{DLSSV05} where a spectral triple on $SU_{2,q}$ is developed.  However, the conventions differ slightly. The coordinate algebra $\mathbb{C}[SU_{2,q}]$ is realized as a sum over the highest weights of the tensor products $M_{\lambda}^* \otimes M_{\lambda}$ where $M_{\lambda}$ is an irreducible module of the highest weight $\lambda$. In the reference \cite{DLSSV05} the duals $M_{\lambda}^*$ are identified with $M_{\lambda}$, which we shall not do here. This leads to a different but isomorphic left representation of $U_q(\mathfrak{su}_2)$. I have decided to follow the usual classical geometric conventions since the identification of the dual modules with the modules is very specific to the case of $\mathfrak{su}_2$ and does not hold in general in the theory of semisimple Lie groups. In addition, with the conventions applied in this work, the spectral triple fits into the general theory of spectral triples on compact quantum groups in \cite{NT10}. \3

I wish to thank the referee for valuable comments which led to significant improvements to the manuscript. \3

\noindent \textbf{Notation.} The parameter $q$ will denote a real number in $(0,1)$. The $q$-integers will be used in the representation theory: 
\begin{eqnarray}\label{qinteger}
[n]_q = \frac{q^n - q^{-n}}{q-q^{-1}}. 
\end{eqnarray}
For coproducts in Hopf algebras,  Sweedler's notation is used: $\triangle(x) = x' \otimes x''$. Whenever the symbol $\updownarrows$ appears it should be understood that the formula is true for both indices $\uparrow$ and $\downarrow$. I shall exploit the notation and refer to weights in the irreducible representations of $\mathfrak{su}_2$ simply by the eigenvalues of the Cartan generator instead of an element in the dual Cartan subalgebra. The conventions where the highest weight runs over $\frac{1}{2} \mathbb{N}_0$ will be used. 

\section{Quantum Group Preliminaries}

\noindent \textbf{1.1.} The quantum Drinfeld-Jimbo algebra $U_q(\mathfrak{su}_2)$ is the complex polynomial algebra generated by $e,f,k,k^{-1}$ subject to the relations
\begin{eqnarray*}
ke = q ek, \5 kf = q^{-1}fe, \5 [e,f] = \frac{k^2 - k^{-2}}{q - q^{-1}}, \5 k k^{-1} = k^{-1} k = 1.
\end{eqnarray*}
$U_q(\mathfrak{su}_2)$ is a Hopf algebra with the coproduct 
\begin{eqnarray*}
\triangle(k) = k \otimes k,\5 \triangle(e) = e \otimes k + k^{-1} \otimes e, \5 \triangle(f) = f \otimes k + k^{-1} \otimes f,
\end{eqnarray*}
the antipode 
\begin{eqnarray*}
S(k) = k^{-1}, \5 S(e) = -qe,\5 S(f) = -q^{-1} f
\end{eqnarray*}
and the counit defined by $\epsilon(k) = 1$, $\epsilon(e) = \epsilon(f) = 0$. We also make $U_q(\mathfrak{su}_2)$ a $*$-algebra by setting $e^* = f$, $f^* = e$ and $(k^{\pm 1})^* = k^{\pm 1}$. 

The representation theory of $U_q(\mathfrak{su}_2)$ is parallel to that of $\mathfrak{su}_2$. The irreducible finite dimensional complex representations are parametrized by the highest weight of $\mathfrak{su}_2$. We denote by $M_{\lambda}$ the $U_q(\mathfrak{su}_2)$-module of highest weight $\lambda$. The dimension of $M_{\lambda}$ is equal to $(2\lambda+1)$ and the basis vectors $u_{\lambda m}$, $-\lambda  \leq m \leq \lambda $ are chosen so that the representation takes the ladder operator form
\begin{eqnarray}\label{ladder}
\varrho_{\lambda }(e)u_{\lambda m}  &=& \sqrt{[\lambda -m]_q[\lambda +m+1]_q} u_{\lambda ,m+1} , \\
\varrho_{\lambda }(f)u_{\lambda m}  &=& \sqrt{[\lambda -m+1]_q[\lambda +m]_q} u_{\lambda ,m-1} ,\nonumber \\
\varrho_{\lambda }(k) u_{\lambda m}  &=& q^{m} u_{\lambda m} .\nonumber
\end{eqnarray}
The $q$-integers $[n]_q$ are defined as in \eqref{qinteger}. Each $\varrho_{\lambda}: U_q(\mathfrak{su}_2) \rightarrow \mathcal{B}(M_{\lambda})$ is surjective but never injective. The tensor products of the irreducible representations compose into a sum of irreducible components, as in the case of $\mathfrak{su}_2$:
\begin{eqnarray}\label{tensorproduct}
M_{\lambda} \otimes M_{\lambda'} = \bigoplus_{\mu = |\lambda - \lambda'|}^{\lambda + \lambda '} M_{\mu}. 
\end{eqnarray}
The coproduct is applied for the action of $U_q(\mathfrak{su_2})$ on the tensor product. Since the coproduct is noncocommutative the symmetric group does not act on the tensor products. This is where the braiding comes in: the Artin's braid group can be given an action on the tensor product. The Clebsch-Gordan coefficients have $q$-deformations which can be used to put the tensor product module to the ladder operator form. \3

\noindent \textbf{1.2.} The enveloping algebra $U_q(\mathfrak{su}_2)$ has a Hopf dual algebra which is spanned as a vector space by the matrix elements of the irreducible representations of $U_q(\mathfrak{su}_2)$:
\begin{eqnarray*}
\mathbb{C}[SU_{2,q}] = \bigoplus_{\lambda \in \frac{1}{2} \mathbb{N}_0} M^*_{\lambda} \otimes M_{\lambda}.
\end{eqnarray*}
For each highest weight $\lambda$ we denote by $u^*_{\lambda m}$, $-\lambda \leq m \leq \lambda$ the dual basis of $M_{\lambda}$. The duality is fixed by $u_{\lambda m}^*(u_{\lambda n}) = \delta_{mn}$. Define
\begin{eqnarray*}
t^{\lambda}_{mn} = u^*_{\lambda m} \otimes u_{ \lambda n}.
\end{eqnarray*}
The dual pairing of $U_q(\mathfrak{su}_2)$ with $\mathbb{C}[SU_{2,q}]$ is defined by  
\begin{eqnarray}\label{pair}
t^{\lambda}_{mn}(x) = u^*_{\lambda m}(\varrho_{\lambda}(x) u_{\lambda n}).
\end{eqnarray}
The dual space $M_{\lambda}^*$ becomes a $U_q(\mathfrak{su}_2)$ module under the dual representation which is defined by: 
\begin{eqnarray*}
\varrho^*_{\lambda}(x) = (\varrho_\lambda(S(x)))^{t} 
\end{eqnarray*}
for all $x \in U_q(\mathfrak{su}_2)$. 

The product in the Hopf dual $\mathbb{C}[SU_{2,q}]$ is defined by requiring that 
\begin{eqnarray}\label{product}
t^{\lambda}_{mn} t^{\lambda'}_{m'n'}(x) = t^{\lambda}_{mn}(x') t^{\lambda'}_{m'n'}(x'')
\end{eqnarray}
for all $x \in U_q(\mathfrak{su}_2)$ and for all possible weight parameters $\lambda,m,n$. \3

\noindent \textbf{Lemma 1.} The assignment \eqref{product} makes $\mathbb{C}[SU_{2,q}]$ a unital complex algebra with the unit $t^0_{00}$ and the product is given by 
\begin{eqnarray}\label{product2}
t^{\lambda}_{mn} t^{\lambda'}_{m'n'} = \sum_{\mu = |\lambda-\lambda'|}^{\lambda +\lambda'}C_q  \begin{pmatrix} \lambda & \lambda' & \mu \\ m & m' & m + m' \end{pmatrix} C_q  \begin{pmatrix} \lambda & \lambda' & \mu \\ n & n' & n+n' \end{pmatrix} t^{\mu}_{m+m',n+n'}. 
\end{eqnarray}
where $C_q(\cdot)$ denote the orthogonal Clebsch-Gordan matrices for the representations $(\varrho_{\lambda}, M_{\lambda})$. \3

\noindent Proof. The Clebsch-Gordan matrices $C_q$ can be applied to write the tensor product components in terms of weight vectors of some modules that appear in the tensor product decomposition. We choose a basis so that $C_q$ are orthogonal. The right side of \eqref{product} can be manipulated by 
\begin{eqnarray*}
&& u_{\lambda m}^* \otimes u_{\lambda n}^* ((\varrho_{\lambda}(x') \otimes \varrho_{\lambda}(x'')) u_{\lambda  m'} \otimes u_{\lambda n'})  \\ &=& (C_q u_{\lambda  m} \otimes u_{\lambda n})^* C_q(\varrho_{\lambda }(x') \otimes \varrho_{\lambda }(x''))C_q^{-1} (C_q u_{\lambda m'} \otimes u_{\lambda n'})
\end{eqnarray*}
Now $C_q(\varrho_{\lambda }(x') \otimes \varrho_{\lambda }(x''))C_q^{-1}$ is the representation of $x$ on the new basis where the irreducible components are put in the ladder operator form by acting on $C_q$. When written out, we get \eqref{product2}. It is elementary to check that $t^0_{00}$ is the unit. \5 $\square$\3

 The product in $U_q(\mathfrak{su}_2)$ is obviously compatible with the coproduct in $\mathbb{C}[SU_{2,q}]$ defined by 
\begin{eqnarray*}
\triangle(t^{\lambda}_{mn}) = \sum_k t^{\lambda}_{mk} \otimes t^{\lambda}_{kn},
\end{eqnarray*}
The antipode and the counit in the Hopf dual $\mathbb{C}[SU_{2,q}]$ are determined by 
\begin{eqnarray*}
t(S(x)) = (S(t))(x),\5 t(1) = \epsilon(t), \5 1(x) = t^0_{00}(x) = \epsilon(x)
\end{eqnarray*}
for all $t \in \mathbb{C}[SU_{2,q}]$ and $x \in U_q(\mathfrak{su}_2)$. Note that same symbols are used for the Hopf-algebra structure maps in $U_q(\mathfrak{su}_2)$ and $\mathbb{C}[SU_{2,q}]$ but this should not lead to confusion. The Hopf dual $\mathbb{C}[SU_{2,q}]$ is equipped with the $*$-structure: $t^*(x) = \overline{t(S(x)^*)}$. 

As a $*$-algebra $\mathbb{C}[SU_{2,q}]$ is generated by the elements $\alpha$ and $\beta$ which, together with their involutions, are given by 
\begin{eqnarray*}
\alpha = t^{\frac{1}{2}}_{\frac{1}{2},\frac{1}{2}}, \5 \beta = t^{\frac{1}{2}}_{\frac{1}{2},-\frac{1}{2}}, \5 \alpha^* = t^{\frac{1}{2}}_{-\frac{1}{2},-\frac{1}{2}}, \5 \beta^* = - \frac{1}{q} t^{\frac{1}{2}}_{-\frac{1}{2},\frac{1}{2}}.
\end{eqnarray*}
In terms of generators and relations, the algebraic structure of $\mathbb{C}[SU_{2,q}]$ is determined by 
\begin{eqnarray*} 
&& \beta \alpha = q\alpha \beta  , \5 \beta^* \alpha = q \alpha \beta^* , \5 \beta \beta^* = \beta^* \beta \\
&&\alpha \alpha^* + \beta \beta^* = 1, \5 \alpha^* \alpha + q^2 \beta^* \beta = 1
\end{eqnarray*}
which is elementary to check by applying the Clebsch-Gordan matrices which are given in \cite{KS97}. The dual pairing applied to the generators $e,f,k,k^{-1}$ of $U_q(\mathfrak{su}_2)$ and to the elements $\alpha, \alpha^*, \beta, \beta^*$ gives
\begin{eqnarray*}
\alpha(k^{\pm 1}) = q^{\pm \frac{1}{2}}, \5 \alpha^*(k^{\pm 1}) = q^{\mp \frac{1}{2}}, \5 \beta(e) = 1, \5 \beta^*(f) = - \frac{1}{q}.
\end{eqnarray*}
and zero in the remaining cases of pairings between these elements. Explicit formulas for the Hopf algebra maps can be solved by dualizing the Hopf structure in $U_q(\mathfrak{su}_2)$: 
\begin{eqnarray*}
&&\triangle(\alpha) = \alpha \otimes \alpha - q \beta \otimes \beta^*, \5 \triangle(\beta) = \beta \otimes \alpha^* + \alpha \otimes \beta, \\ 
&&S(\alpha) = \alpha^*, \5 S(\beta) = -q\beta, \5 S(\beta^*) = - \frac{1}{q} \beta^*, \5 S(\alpha^*) = \alpha, \\
&&\epsilon(\alpha) = 1, \5 \epsilon(\beta) = 0. 
\end{eqnarray*}
We make $U_q(\mathfrak{su}_2)$ act on  $\mathbb{C}[SU_{2,q}]$ under the right and the left regular representations:  
\begin{eqnarray}\label{lr}
\partial(x)(t^{\lambda}_{mn}) = u^*_{\lambda m} \otimes \varrho_{\lambda}(x)u_{\lambda n}, \5 l(x)(t^{\lambda}_{mn}) = \varrho^*_{\lambda}(\vartheta (x)) u^*_{\lambda m}  \otimes  u_{\lambda n} 
\end{eqnarray}
for all $x \in U_q(\mathfrak{su}_2)$ and $t^{\lambda}_{mn} \in \mathbb{C}[SU_{2,q}]$. The left representation is equipped with the following algebra automorphism
\begin{eqnarray*}
\vartheta(k^{\pm 1}) = k^{\mp 1}, \5 \vartheta(e) = - f, \5 \vartheta(f) = -e.
\end{eqnarray*}
This automorphism restores the commutativity of the coproduct and the antipode in the sense that  
\begin{eqnarray}\label{theta}
\triangle(S(\vartheta(x)) = S(\vartheta(x')) \otimes S(\vartheta(x'')) 
\end{eqnarray}
for all $x \in U_q(\mathfrak{su}_2)$. When applied on the generators of $\mathbb{C}[SU_{2,q}]$ the right regular representation is given by 
\begin{center}
\begin{tabular}{l l l}
$\partial(e) \alpha = 0 $& $\partial(f) \alpha =  \beta$ & $\partial(k^{\pm 1}) \alpha = q^{\pm \frac{1}{2}} \alpha$ \\
$\partial(e) \beta = \alpha$ & $\partial(f) \beta = 0$ & $\partial(k^{\pm 1}) \beta = q^{\mp \frac{1}{2}} \beta$ \\
$\partial(e) \alpha^* = -q\beta^*$ & $\partial(f) \alpha^* = 0$ & $\partial(k^{\pm 1}) \alpha^* = q^{\mp \frac{1}{2}} \alpha^*$ \\
$\partial(e) \beta^* = 0$ & $\partial(f) \beta^* = - \frac{1}{q} \alpha^*$ & $\partial(k^{\pm 1}) \beta^* = q^{\pm \frac{1}{2}} \beta^*$. 
\end{tabular}
\end{center}
For the left representation one needs to apply the antipode and work with the dual representation. A straightforward computations gives: 
\begin{center}
\begin{tabular}{l l l}
$l(e) \alpha = 0$ & $l(f) \alpha = - q^2 \beta^*$ & $l(k^{\pm 1}) \alpha = q^{\pm \frac{1}{2}} \alpha$ \\
$l(e) \beta = 0$ & $l(f) \beta = q \alpha^*$ & $l(k^{\pm 1}) \beta = q^{\pm \frac{1}{2}} \beta$ \\
$l(e) \alpha^* = \frac{1}{q} \beta$ & $l(f) \alpha^* = 0$ & $l(k^{\pm 1}) \alpha^* = q^{\mp \frac{1}{2}} \alpha^*$ \\
$l(e) \beta^* = - \frac{1}{q^2} \alpha$ & $l(f) \beta^* = 0$ & $l(k^{\pm 1}) \beta^* = q^{ \mp \frac{1}{2}} \beta^*$.
\end{tabular}
\end{center}

\noindent \textbf{1.3.} The Haar state in the algebra $\mathbb{C}[SU_{2,q}]$ is a linear functional $h: \mathbb{C}[SU_{2,q}] \rightarrow \mathbb{C}$ that is fixed by the relations $h(1) = 1$ and $h(t^{\lambda}_{mn}) = 0$ for all ${\lambda} > 0$. The Haar state provides a Hilbert space completion for $\mathbb{C}[SU_{2,q}]$, which will be denoted by $L^2(SU_{2,q})$, and the GNS construction defines a representation $\pi_h$ of $\mathbb{C}[SU_{2,q}]$ on $L^2(SU_{2,q})$. The basis vectors are mutually orthogonal \cite{KS97}: 
\begin{eqnarray*}
h((t^{\lambda}_{mn})^* t^{\lambda'}_{m'n'}) = \frac{q^{-2m}}{[2\lambda+1]_q} \delta_{\lambda \lambda'} \delta_{mm'} \delta_{nn'}.  
\end{eqnarray*}
Let us denote by $\eta: \mathbb{C}[SU_{2,q}] \rightarrow L^2(SU_{2,q})$ the natural inclusion. Now $\eta(t^{\lambda}_{mn})$ gives a basis of $L^2(SU_{2,q})$. The representation $\pi_h : \mathbb{C}[SU_{2,q}] \rightarrow L^2(SU_{2,q})$ is given by 
\begin{eqnarray*}
\pi_h(t^{\lambda}_{mn}) \eta(t^{\lambda'}_{m'n'})   = \eta(t^{\lambda}_{mn} t^{\lambda'}_{m'n'}),
\end{eqnarray*}
recall the product rule \eqref{product2}. An orthonormal basis of the Hilbert space is given by 
\begin{eqnarray*}
|\lambda m n \ra = q^m \sqrt{[2 \lambda + 1]_q} t^{\lambda}_{mn}. 
\end{eqnarray*}
Explicit formulas for the action of the generators of $\mathbb{C}[SU_{2,q}]$ in this basis is computed in the reference \cite{DLSSV05} (the only notational difference is that they use the symbol $l$ for the highest weight). 

The right and the left representations of $U_q(\mathfrak{su}_2)$ on $\mathbb{C}[SU_{2,q}]$ can be extended on $L^2(SU_{2,q})$ by letting the quantum group act on the basis according to 
\begin{eqnarray*}
\partial(x) \eta(t^{\lambda}_{mn}) = \eta(\partial(x)t^{\lambda}_{mn}),\5 l(x) \eta(t^{\lambda}_{mn}) = \eta(l(x)t^{\lambda}_{mn}) 
\end{eqnarray*}
for all $x \in U_q(\mathfrak{su}_2)$ and $t^{\lambda}_{mn} \in \mathbb{C}[SU_{2,q}]$. \3

\noindent \textbf{Proposition 1.} The GNS representation of $\mathbb{C}[SU_{2,q}]$ is equivariant under the right and the left regular representations of $U_q(\mathfrak{su}_2)$ in the sense that 
\begin{eqnarray*}
\partial(x)(\pi_h(t^{\lambda}_{mn}) \eta(t^{\lambda'}_{m'n'})) &=& \pi_h(\partial(x') t^{\lambda}_{mn}) \eta(\partial(x'')t^{\lambda'}_{m'n'}) \\
l(x)(\pi_h(t^{\lambda}_{mn}) \eta(t^{\lambda'}_{m'n'})) &=& \pi_h(l(x') t^{\lambda}_{mn}) \eta(l(x'')t^{\lambda'}_{m'n'})
\end{eqnarray*}
for all $x \in U_q(\mathfrak{su}_2)$. \3

\noindent Proof. The pairing of $t \in  \mathbb{C}[SU_{2,q}]$ and $x \in U_q(\mathfrak{su}_2)$ takes the matrix element $t$ of the right representation $\partial(x)$ which is evident from the definitions \eqref{pair}, \eqref{lr}. Two right representations are equal if and only if all their matrix components are equal. Therefore the right equivariance follows from the definition of the product \eqref{product}. The left equivariance is a consequence of
\begin{eqnarray*}
t^{\lambda}_{mn} t^{\lambda'}_{m'n'} (S\vartheta(x)) &=& t^{\lambda}_{mn}((S\vartheta(x))') t^{\lambda'}_{m'n'} ((S\vartheta(x))'') \\ 
&=&  t^{\lambda}_{mn}((S\vartheta(x'))) t^{\lambda'}_{m'n'} ((S\vartheta(x''))) 
\end{eqnarray*}
where the first equality follows from the definition of the product and the second follows from \eqref{theta}. \5 $\square$ \3

\noindent \textbf{1.4.} The next step is to define a Hilbert space of spinors. In the geometric model over $SU_2$ the complexified Clifford algebra over the Lie algebra $\mathfrak{su}_2$ is semisimple and isomorphic to two copies of $\mathcal{B}(\mathbb{C}^2)$. The spinor module is an irreducible representation of the Clifford algebra and therefore two dimensional. So it makes sense to define the spinor module in the quantum group model to be $M_{\frac{1}{2}}$. Denote by $e_{\pm}$ the weight $\pm \frac{1}{2}$ basis vectors of $M_{\frac{1}{2}}$.  We construct the Hilbert space of spinors over $SU_{2,q}$ by tensoring $L^2(SU_{2,q})$ with $M_{\frac{1}{2}}$. Define 
\begin{eqnarray*}
\hil = L^2(SU_{2,q}) \otimes M_{\frac{1}{2}}. 
\end{eqnarray*}
The right and the left regular representations are defined by
\begin{eqnarray*}
\partial'(x) = \partial(x') \otimes \varrho_{\frac{1}{2}}(x''), \5 l'(x) = l(x) \otimes \iota. 
\end{eqnarray*}
We decompose $\hil$ to the irreducible components with respect to these actions. Under the right representation, one has
\begin{eqnarray*}
M_0 \otimes M_{\frac{1}{2}} = M_{\frac{1}{2}} \5 \text{and}\5 M_{\lambda} \otimes M_{\frac{1}{2}} = M_{\lambda - \frac{1}{2}} \oplus M_{\lambda + \frac{1}{2}} \5 \text{if} \5 \lambda > 0. 
\end{eqnarray*}
To extract the irreducible components in the tensor product we need to change the basis. We follow \cite{DLSSV05} and define a new basis in which the tensor product components take the ladder operator form. Define $k^{\pm} = k \pm \frac{1}{2}$ for all $k \in \frac{1}{2} \mathbb{Z}$. 

For $j \in \frac{1}{2}\mathbb{N}$ we define the vectors 
\begin{eqnarray*}
|j \mu \downarrow \ra = C_{j \mu} (u_{j^- \mu^+} \otimes e_-) + S_{j \mu} (u_{j^- \mu^-} \otimes e_+) 
\end{eqnarray*}
where $\mu \in \{-j, \ldots,j-1, j\}$. For $j \in \frac{1}{2} \mathbb{N}_0$ we define
\begin{eqnarray*} 
|j \mu \uparrow \ra = -S_{j+1, \mu} (u_{j^+ \mu^+} \otimes e_-)+ C_{j+1, \mu} (u_{j^+ \mu^-} \otimes e_+)
\end{eqnarray*}
where $\mu \in \{-j, \ldots,j-1, j\}$. The coefficients are given by 
\begin{eqnarray*}
C_{j \mu} = q^{\frac{-(j + \mu)}{2}} \frac{[j- \mu]_q^{\frac{1}{2}}}{[2j]_q^{\frac{1}{2}}}, \5 S_{j \mu} = q^{\frac{j  - \mu}{2}} \frac{[j+ \mu]_q^{\frac{1}{2}}}{[2j]_q^{\frac{1}{2}}}.
\end{eqnarray*}
The basis of the Hilbert space $\hil$ will be labeled by
\begin{eqnarray*}
|j  m \mu \downarrow \ra &=& q^m \sqrt{[2j^- + 1]_q} (\eta \otimes \iota) (u^*_{j^-, m} \otimes |j \mu \downarrow \ra) \\
&=& q^m \sqrt{[2j^- + 1]_q} (C_{j \mu} \eta(t^{j^-}_{m, \mu^+}) \otimes e_- + S_{j \mu} \eta(t^{j^-}_{m, \mu^-}) \otimes e_+), \\ 
|j  m \mu \uparrow \ra &=& q^m \sqrt{[2j^+ + 1]_q} (\eta \otimes \iota)(u^*_{j^+, m} \otimes |j \mu \uparrow \ra) \\ &=&  q^m \sqrt{[2j^+ + 1]_q}(-S_{j+1, \mu} \eta(t^{j^+}_{m,\mu^+}) \otimes e_- + C_{j+1,\mu} \eta(t^{j^+}_{m,\mu^-}) \otimes e_+)  
\end{eqnarray*}
where $m$ runs over the usual weight space parameters: $-j^{-} \leq m \leq j^{-}$ in $|j  m \mu \downarrow \ra$ and $-j^+ \leq m \leq j^+$ in  $|j  m \mu \uparrow \ra$. Moreover, $-j \leq \mu \leq j$ in both cases. 

The basis is orthonormal which is a consequence of the property $C^2_{j \mu} + S^2_{j\mu} = 1$ for the Clebsch-Gordan coefficients. In this basis the right regular representation takes the ladder operator form 
\begin{eqnarray*}
\partial'(x) |j m \mu \downarrow \ra &=& \sum_{\nu} \varrho_{j^-}(x)_{\nu \mu} |j m \nu \downarrow \ra, \\
\partial'(x) |j m \mu \uparrow \ra &=& \sum_{\nu} \varrho_{j^+}(x)_{\nu \mu} |j m \nu \uparrow \ra
\end{eqnarray*}
where $\varrho_j(x)_{\nu \mu}$ are the matrix coefficients of the operation $\varrho_j(x)$,  recall \eqref{ladder}. 

On the Hilbert space $\hil$ we use the representation of $\mathbb{C}[SU_{2,q}]$ defined by $\pi = \pi_h \otimes \iota$. \3

\noindent \textbf{1.5.} The construction of a Dirac operator in \cite{DLSSV05}, or in \cite{NT10}, is based on the invariance under the right and the left representations of $U_q(\mathfrak{su}_2)$. Recall that all the irreducible representations of $U_q(\mathfrak{su}_2)$ are surjective maps $U_q(\mathfrak{su}_2) \rightarrow \mathcal{B}(M_{\lambda})$. So, in order to make the Dirac operator commute with both actions, one needs to require that it acts constantly on each irreducible component. Thus $D$ is diagonal in the basis $|j m \mu \updownarrows \ra$ and the spectrum can only depend on the parameter $j$. We shall apply the classical Dirac spectrum. Explicitly, this operator is an unbounded densely defined self-adjoint operator which acts on the basis of $\hil$ by 
\begin{eqnarray*}
D |j m \mu \uparrow \ra = (2j + \frac{3}{2})|j  m \mu \uparrow \ra, \5 D |j m \mu \downarrow \ra = -(2j+\frac{1}{2}) |j m \mu \downarrow \ra. 
\end{eqnarray*}
The multiplicities of these eigenvalues are $(2j+1)(2j+2)$ and $2j(2j+1)$. In the classical limit $q \rightarrow 1$ this model coincides with the Dirac operator associated to the bi-invariant metric and Levi-Civita connection on $SU_2$, \cite{Sle85}. 

A more general model for spinors and Dirac operators on compact quantum groups is developed in \cite{NT10}. In this case the quantum group Dirac operator is a unitary conjugate of the classical Dirac operator and so the classical spectrum is automatically preserved. In addition this operator has the right and the left symmetries under the quantum group representations on the spinor space. On the technical level one needs to apply certain Drinfeld's twist to write down such an operator. The case of $\mathfrak{su}_2$ is special in the theory of semisimple Lie algebras since the irreducible components in the tensor product decomposition \eqref{tensorproduct} appear always with the multiplicity one. This means that the highest weight vectors in each tensor product $M_{\lambda} \otimes M_{\frac{1}{2}}$ are uniquely determined after the irreducible representations are fixed. Therefore, the decomposition into irreducible components and consequently the decomposition into the Dirac eigenspaces is unique. It would be reasonable to expect that the model of \cite{NT10} would coincide with the construction above. Indeed, with the conventions applied above this happens. In \cite{NT10} it was noted that after some algebraic manipulation  (and after a normalization), the Dirac operator $\eth$ satisfies
\begin{eqnarray*}
q^{- \eth} = q^{\frac{3}{2}}\begin{pmatrix} \partial(k^2-q^{-1}(q-q^{-1})^2fe) & q^{- \frac{1}{2}}(q-q^{-1})\partial(fk^{-1}) \\
q^{-\frac{1}{2}}(q-q^{-1})\partial(k^{-1}e) & \partial(k^{-2}) \end{pmatrix}
\end{eqnarray*}
It is straightforward to check that 
\begin{eqnarray*}
q^{- \eth} |j m \mu \uparrow \ra = q^{-(2j + \frac{3}{2})}|j  m \mu \uparrow \ra, \5 q^{- \eth} |j m \mu \downarrow \ra = q^{2j+\frac{1}{2}} |j m \mu \downarrow \ra.
\end{eqnarray*}
So, the Dirac operators $D$ and $\eth$ are the same densely defined self-adjoint operators on the spinor module since their eigenspace decomposition and the eigenvalues match. In particular, $(\mathbb{C}[SU_{2,q}], \pi, \hil, D)$ is a spectral triple and the commutators $[D, \pi(t)]$ extend to bounded operators on $\hil$ for all $t \in \mathbb{C}[SU_{2,q}]$. 

\section{Spinors on Quantum Weighted Projective Spaces} 
 
\noindent \textbf{2.1.} The algebras of quantum projective planes are coinvariant subalgebras in $\mathbb{C}[SU_{2,q}]$ under the coaction of the polynomial algebra $B:=\mathbb{C}[u,u^{-1}]$. There is a Hopf $*$-algebra structure in $B$ so that $u$ is a grouplike element and unitary. For each pair of positive coprime integers $k,l$ define the coaction $\theta : \mathbb{C}[SU_{2,q}] \rightarrow \mathbb{C}[SU_{2,q}] \otimes B$ by
\begin{eqnarray*}
\alpha \mapsto \alpha \otimes u^{-k}, \5 \beta \mapsto \beta \otimes u^{l}. 
\end{eqnarray*}
This coaction is exactly the same that was used in \cite{BF12}: the algebra elements denoted by $\alpha$ and $\beta$ in the reference \cite{BF12} are the elements for which we used the symbols $\alpha^*$ and $q \beta^*$ in 1.2. The coinvariant subalgebras are denoted by $\mathbb{C}[\mathbb{WP}_{k,l,q}]$.  It is elementary to check that $\mathbb{C}[\mathbb{WP}_{k,l,q}]$ is the $*$-subalgebra of $\mathbb{C}[SU_{2,q}]$ generated by  
\begin{eqnarray*}
a =  \beta \beta^* \5 \text{and} \5 b = \beta^k \alpha^l.  
\end{eqnarray*}
In terms of generators and relations, the algebra $\mathbb{C}[\mathbb{WP}_{k,l,q}]$ can be identified with the complex polynomial $*$-algebra generated by the symbols $a,b$  which are subject to the relations \cite{BF12}: 
\begin{eqnarray*}
&&a^* = a, \5 ab^* = q^{-2l} b^* a, \5 b^* b = q^{2kl} a^k \prod_{m = 0}^{l-1} (1-q^{2(m+1)}a), \\
&&b b^* = a^k \prod_{m = 1}^l(1-q^{-2(m-1)} a) \nonumber. 
\end{eqnarray*}
The coinvariant algebra $\mathbb{C}[\mathbb{WP}_{k,l,q}]$ is not equipped with a coalgebra structure.\3
 
\noindent \textbf{2.2.} The generators of the algebra $\mathbb{C}[SU_{2,q}]$ are homogeneous elements under the coaction $\theta$ and so there is the direct sum decomposition
\begin{eqnarray}\label{module}
\mathbb{C}[SU_{2,q}] = \bigoplus_{i \in \mathbb{Z}} \mathbb{C}[SU_{2,q}]^{(i)}
\end{eqnarray}
into the subspaces spanned by the homogeneous vectors of order $i$. The components are orthogonal with respect to the Haar state of 1.3: the vectors $t^{\lambda}_{mn}$ define an orthogonal basis for $\mathbb{C}[SU_{2,q}]$ and these vectors are homogeneous.\3

\noindent \textbf{Proposition 2.} The basis vectors $t^{\lambda}_{mn}$ of $\mathbb{C}[SU_{2,q}]$ are homogeneous under the coaction and the coaction takes the form 
\begin{eqnarray*}
\theta(t^{\lambda}_{mn}) = t^{\lambda}_{mn} \otimes u^{-(m+n)k + (m-n)l}
\end{eqnarray*}
for all highest weights $\lambda$ and $-\lambda \leq m,n \leq \lambda$. \3
 
\noindent Proof. The proof follows easily from the formulas in 2.1.  \5 $\square$ \3

In the notation \eqref{module}, $\mathbb{C}[\mathbb{WP}_{k,l,q}]$ is the component of homogeneous order $0$, and the product in $\mathbb{C}[SU_{2,q}]$ determines a $*$-representation of $\mathbb{C}[\mathbb{WP}_{k,l,q}]$ on every $\mathbb{C}[SU_{2,q}]^{(i)}$. So, \eqref{module} has a structure of a direct sum of $\mathbb{C}[\mathbb{WP}_{k,l,q}]$-modules.\3
 
\noindent \textbf{2.3.} The goal in the following is to understand the structure of \eqref{module} as a representation space for $\mathbb{C}[\mathbb{WP}_{k,l,q}]$. Let us discuss the topological case first. Consider the case where $\mathbb{T}$ acts on $SU_2$ from the right with trivial isotropy so that the orbit space is homeomorphic to the 2-sphere $S^2$. Let us realize this action by $(\sigma, g) \mapsto g \sigma^{-1}$ for all $g \in SU_2$ and $\sigma \in \mathbb{T}$. Suppose that for each highest weight of $SU_2$ we have a fixed irreducible finite dimensional $SU_2$-representation. Then the space $\mathbb{C}[SU_2]$ of matrix elements of the finite dimensional irreducible $SU_2$ representations has the direct sum decomposition 
\begin{eqnarray*}
\mathbb{C}[SU_2] = \bigoplus_{i \in \mathbb{Z}} \mathbb{C}[SU_2]^{(i)}  
\end{eqnarray*}
so that the components $\mathbb{C}[SU_2]^{(i)}$ are the homogeneous spaces under the coaction determined by the $\mathbb{T}$ action on $SU_2$. More precisely, $\mathbb{C}[SU_2]^{(i)}$ is the homogeneous subspace of order $i$ under the coaction of $B$ which is defined on the generators by 
\begin{eqnarray*}
t^{\frac{1}{2}}_{\frac{1}{2}, \frac{1}{2}} \mapsto t^{\frac{1}{2}}_{\frac{1}{2}, \frac{1}{2}} \otimes u^{-1} \5 t^{\frac{1}{2}}_{\frac{1}{2}, -\frac{1}{2}} \mapsto t^{\frac{1}{2}}_{\frac{1}{2}, -\frac{1}{2}} \otimes u.
\end{eqnarray*}
These generators are matrix elements of the irreducible $2$-dimensional representation of $SU_2$ and the parameters $\pm \frac{1}{2}$ correspond to the weights of this module, exactly as in the case of the quantum group in 1.2. The completion of $\mathbb{C}[SU_2]^{(0)}$ in the uniform norm is isomorphic to the $C^*$-algebra of continuous functions on $S^2$, which we denote by $C(SU_2)^{(0)}$. The completions $\overline{\mathbb{C}[SU_2]^{(i)}}$ in the uniform norm are $C^*$-modules for $C(SU_2)^{(0)}$ for all $i \in \mathbb{Z}$. Under the identification of $C(SU_2)^{(0)}$ with $C(S^2)$, the components $\overline{\mathbb{C}[SU_2]^{(i)}}$ correspond to the projective $C(S^2)$-modules of continuous sections of the complex line bundles on $S^2$. One can use this decomposition to give a full analytic realization for the $C^*$-algebraic K-theory group $K_0(C(S^2)) = \mathbb{Z} \oplus \mathbb{Z}$. The same holds for the teardrop orbifolds which are homeomorphic to $S^2$ and so the (nonequivariant) K-theory groups are the same. We shall study the noncommutative quantum teardrops, i.e. the cases with $k = 1$ and $l \in \mathbb{N}$, and show that the components of \eqref{module} give rise to all projective modules for a $C^*$-algebra completion of $\mathbb{C}[\mathbb{WP}_{1,l,q}]$ and we show that these modules can be applied to give a full analytical realization for the even $C^*$-algebraic K-theory group of the quantum teardrops.

The  irreducible infinite dimensional $*$-representations of the quantum group $\mathbb{C}[SU_{2,q}]$ are parametrized by $t \in \mathbb{T}$ and given by   
\begin{eqnarray*}
\varpi_t(\alpha) e_n = \sqrt{1 - q^{2(n+1)}} e_{n+1}, \5 \varpi_t(\beta) e_n = tq^{n} e_n
\end{eqnarray*}
where $\{e_n: n \in \mathbb{N}_0\}$ is the basis of $l^2(\mathbb{N}_0)$. We can merge these to a single representation on the space $L^2(\mathbb{T}) \otimes l^2(\mathbb{N}_0)$, which can be identified with $l^2(\mathbb{Z}) \otimes l^2(\mathbb{N}_0)$ through the Fourier transform. This gives the faithful $*$-representation: 
\begin{eqnarray*}
\varpi(\alpha) = \iota \otimes \varpi_1(\alpha), \5 \varpi(\beta) = \mathcal{U} \otimes \varpi_1(\beta)
\end{eqnarray*}
where $\mathcal{U}$ is the unilateral forward shift operator on $l^2(\mathbb{Z})$. We shall identify $l^2(\mathbb{N}_0)$ with  $\bigoplus_{s = 1}^l l^2(\mathbb{N}_0)$ through the linear isomorphism
\begin{eqnarray*}
l^2(\mathbb{N}_0) \rightarrow \bigoplus_{s = 1}^l l^2(\mathbb{N}_0), \5  e_{lp + s - 1} \mapsto e^s_p.
\end{eqnarray*}
where $\{e^s_p: p \in \mathbb{N}_0\}$ is the basis of the $s$'th copy of $l^2(\mathbb{N}_0)$. We also set 
\begin{eqnarray*}
l^2(\mathbb{Z}) \otimes \bigoplus_{s = 1}^l l^2(\mathbb{N}_0) &=& \bigoplus_{m \in \mathbb{Z}} \mathcal{X}_m = \bigoplus_{m \in \mathbb{Z}} \bigoplus_{s = 1}^l \mathcal{X}_m^s
\end{eqnarray*}
so that $\mathcal{X}_m$ is the subspace spanned by $e_{p+m} \otimes e_p^s$ for all $p \geq 0$ and $1 \leq s \leq l$, and $\mathcal{X}_m^s$ is the subspace of $\mathcal{X}_m$ for the fixed $s$.

Let us consider the case of quantum teardrops, $\mathbb{C}[\mathbb{WP}_{1,l,q}]$. For each $m \in \mathbb{Z}$, the representation $\varpi$ of $\mathbb{C}[SU_{2,q}]$ restricts to a direct sum of irreducible faithful $*$-representations for the teardrop $*$-subalgebra, 
\begin{eqnarray*}
\varpi^m = \bigoplus_{s = 1}^l \varpi^m_s: \mathbb{C}[\mathbb{WP}_{1,l,q}] \rightarrow \mathcal{B}(\bigoplus_{s = 1}^l \mathcal{X}_m^s),
\end{eqnarray*} 
which are given by 
\begin{eqnarray*}
\varpi^m_s(a) e_{m+p} \otimes e^s_p &=& q^{2(lp+s-1)} e_{m+p} \otimes e^s_p, \\
\varpi^m_s(b^*) e_{m+p} \otimes e^s_p &=& q^{lp + s - 1} \prod_{r = 1}^{l} \sqrt{1 - q^{2(lp + s - r)}} e_{m + p-1} \otimes e^s_{p-1}, \\
\varpi^m_s(b^*) e_{m} \otimes e^s_0 &=& 0.
\end{eqnarray*}
We identify each $\mathcal{X}_m^s$ with $l^2(\mathbb{N}_0)$ by
\begin{eqnarray}\label{hilbid}
e_{m+p} \otimes e^s_p \mapsto e_p. 
\end{eqnarray}
Now it is evident that, for any choice of $m \in \mathbb{Z}$, $\varpi_s^m: \mathbb{C}[\mathbb{WP}_{1,l,q}] \rightarrow \mathcal{B}(\mathcal{X}_m^s) = \mathcal{B}(l^2(\mathbb{N}_0))$ are representatives of the unitary equivalence classes of all bounded faithful irreducible $*$-representations of $\mathbb{C}[\mathbb{WP}_{1,l,q}]$ which were classified in \cite{BF12}. 

The $C^*$-algebra $C(\mathbb{WP}_{1,l,q})$ is the subalgebra of $\mathcal{B}(\bigoplus_{s = 1}^l l^2(\mathbb{N}_0))$ obtained by taking the completion of $\bigoplus_{s = 1}^l (\varpi^m_s(\mathbb{C}[\mathbb{WP}_{1,l,q}]))$ in the operator norm. Let $\mathcal{K}$ denote the $C^*$-algebra of compact operators on $l^2(\mathbb{N}_0)$. Then there is an identification of the $C^*$-algebras \cite{BF12}:
\begin{eqnarray*}
C(\mathbb{WP}_{1,l,q}) = (\bigoplus_{s = 1}^l \mathcal{K} )^+
\end{eqnarray*}
where $+$ denotes the unitization of an algebra. It follows that the even K-theory of this $C^*$-algebra is subject to an isomorphism 
\begin{eqnarray*}
K_0(C(\mathbb{WP}_{1,l,q})) \simeq \mathbb{Z}^{\oplus l + 1}. 
\end{eqnarray*}
The unitary equivalence classes of the projections in the matrix algebra $M_{\infty}((\bigoplus_{s = 1}^l \mathcal{K} )^+)$ can be represented by the projections $\bigoplus_{j = 1}^l P_{k_j}$
for $k_j \in \mathbb{N}_0$, and 
\begin{eqnarray*}
I_r \oplus (1 - (\bigoplus_{j = 1}^l P_{n_j})) \oplus (\bigoplus_{j = 1}^l P_{m_j})
\end{eqnarray*}
for $r, n_j, m_j \in \mathbb{N}_0$ so that $n_j m_j = 0$ for all $j$ where $I_r$ denotes the identity matrix of $M_r((\bigoplus_{s = 1}^l \mathcal{K} )^+)$ and
\begin{eqnarray*}
P_{n_j} := \sum_{i = 0}^{n_j-1} e^j_{ii} \in \mathcal{K}
\end{eqnarray*}
is the finite rank projection onto the subspace spanned by $e_{0}^j, \ldots, e_{n_j-1}^j$ for all $j \in \{1, \ldots, l\}$, and $P_0 = 0$.  One can apply the group structure in K-theory to give a representative for all the elements in $K_0((\bigoplus_{s = 1}^l \mathcal{K} )^+)$ in terms of these projections. 

The analytic realization of the K-theory groups $K_0(C(\mathbb{WP}_{1,l,q}))$ was studied in \cite{She14} in terms of noncommutative bundles over the base algebra $C(\mathbb{WP}_{1,l,q})$. The quantum lens space, $\mathbb{C}[L_q(l;1,l)]$, is the $*$-subalgebra of $\mathbb{C}[SU_{2,q}]$ generated by $\alpha^l$ and $\beta$. The algebra $\mathbb{C}[L_q(l;1,l)]$ has a structure of a quantum principal $\mathbb{T}$-bundle over $\mathbb{C}[\mathbb{WP}_{1,l,q}]$, \cite{BF12}, however we shall not employ this feature here. Clearly $\mathbb{C}[L_q(l;1,l)]$ is the subalgebra of homogeneous components of order $nl$ in $\mathbb{C}[SU_{2,q}]$ for all $n \in \mathbb{Z}$, and so we can write
\begin{eqnarray*}
\mathbb{C}[L_q(l;1,l)] = \bigoplus_{n \in \mathbb{Z}} \mathbb{C}[SU_{2,q}]^{(nl)} = \bigoplus_{n \in \mathbb{Z}} \mathcal{L}[n]. 
\end{eqnarray*}
In our notation 'the basic quantum line bundles' $\mathcal{L}[n]$, studied in \cite{BF12} and \cite{She14}, are the modules $\mathbb{C}[SU_{2,q}]^{(nl)}$. The representation $\varpi: \mathbb{C}[SU_{2,q}] \rightarrow \mathcal{B}(l^2(\mathbb{Z}) \otimes \bigoplus_{s = 1}^l l^2(\mathbb{N}_0))$ restricts on the subalgebra $\mathbb{C}[L_q(l;1,l)]$ to give the direct sum of $*$-representations 
\begin{eqnarray*}
\varpi = \bigoplus_{s = 1}^l \varpi_s: \mathbb{C}[L_q(l;1,l)] \rightarrow \mathcal{B}(l^2(\mathbb{Z}) \otimes \bigoplus_{s = 1}^l l^2(\mathbb{N}_0)),
\end{eqnarray*}
which are defined by 
\begin{eqnarray}\label{repre}
\varpi_s(\alpha^l) e_{p+m} \otimes e^s_p &=& \prod_{r = 0}^{l-1} \sqrt{1-q^{2(pl + s + r)}} e_{p+m} \otimes e^s_{p+1}, \nonumber \\
\varpi_s(\beta) e_{p+m} \otimes e^s_p &=& q^{lp+s-1} e_{p+ m + 1} \otimes e^s_{p}.  
\end{eqnarray}
for all $m \in \mathbb{Z}$. The components of $\varpi = \bigoplus_{s = 1}^l \varpi_s$ are bounded faithful irreducible $*$-representations of the algebra $\mathbb{C}[L_q(l;1,l)]$. The completion of $\varpi(\mathbb{C}[L_q(l;1,l)])$ in the operator norm is the $C^*$-algebra of the quantum lens space, $C(L_q(l;1,l))$. The completion of $\varpi(\mathcal{L}[n])$ in $C(L_q(l;1,l))$ is denoted by $\overline{\mathcal{L}[n]}$. Each $\overline{\mathcal{L}[n]}$ is a projective module for the $C^*$-algebra, $C(\mathbb{WP}_{1,l,q}) = (\bigoplus_{s = 1}^l \mathcal{K} )^+$, and its K-theory class is represented by the projection, \cite{She14}: 
\begin{eqnarray*}
I_1 \oplus (\bigoplus_{s = 1}^l P_n) \5 \text{if $n \geq 0$} \5 \text{and} \5 1 - \bigoplus_{s = 1}^l P_n \5 \text{if $n < 0$}
\end{eqnarray*}
for all $n \in \mathbb{Z}$.  

Consider the completion of $\varpi(\mathbb{C}[SU_{2,q}]^{(nl+j)})$ in the operator norm, $\overline{\mathbb{C}[SU_{2,q}]^{(nl+j)}}$, as a left $C^*$-module for the algebra $C(\mathbb{WP}_{1,l,q})$ under the operator product. \3

\noindent \textbf{Proposition 3.} The $C(\mathbb{WP}_{1,l,q})$-module $\overline{\mathbb{C}[SU_{2,q}]^{(nl+j)}}$ is isomorphic to 
\begin{eqnarray}\label{module2}
(\bigoplus_{s = 1}^l \mathcal{K}) + \mathbb{C} (\bigoplus_{s = 1}^{l-j} \mathcal{S}^{n} \oplus \bigoplus_{s = l-j+1}^{l} \mathcal{S}^{n+1} )
\end{eqnarray}
as a $C^*$-module for all $l \in \mathbb{N}$, $n \in \mathbb{Z}$ and $j \in \{1, \ldots, l-1\}$. \3
 
We consider \eqref{module2} as a subspace in $\mathcal{B}( \bigoplus_{s = 1}^l l^2(\mathbb{N}_0))$ and the symbol $+$ denotes the vector space sum. $\mathcal{S}$ is the unilateral backward shift on $l^2(\mathbb{N}_0)$. The operator product defines the left module structure. \3

\noindent Proof of Proposition 3. The operator space $\varpi(\mathcal{L}[n])$ restricts to
\begin{eqnarray*}
\varpi(\mathcal{L}[n]): \mathcal{X}_m \rightarrow \mathcal{X}_{m+n}
\end{eqnarray*}
for all $m \in \mathbb{Z}$, recall \eqref{repre}. Identify each $\mathcal{X}_m$ with $l^2(\mathbb{N}_0)$ as in \eqref{hilbid}. The coefficients of $\varpi(\mathcal{L}[n])$ are independent on the parameter $m$ and we can view $\varpi(\mathcal{L}[n])$ as a subspace in $\mathcal{B}(\bigoplus_{s = 1}^l l^2(\mathbb{N}_0))$. In \cite{She14} the completion $\overline{\mathcal{L}[n]}$ in $\mathcal{B}(\bigoplus_{s = 1}^l l^2(\mathbb{N}_0))$ was shown to be the space 
\begin{eqnarray*}
 (\bigoplus_{s = 1}^l \mathcal{K}) + \mathbb{C} (\bigoplus_{s = 1}^{l} \mathcal{S}^n)
\end{eqnarray*}
for all $n \in \mathbb{Z}$. 

The algebra $\mathbb{C}[SU_{2,q}]$ is linearly spanned by $\alpha^{n_1} \beta^{n_2} (\beta^*)^{n_3}$ where $n_1 \in \mathbb{Z}$, $n_2, n_3 \in \mathbb{N}_0$ and $\alpha^{-n} = (\alpha^*)^n$ if $n > 0$, \cite{Wor87b}. So, $\mathbb{C}[SU_{2,q}]^{(ln + j)}$ is the vector space 
\begin{eqnarray}\label{prop3}
(\alpha^*)^j  \mathcal{L}[n] + \alpha^{l-j} \mathcal{L}[n+1]
\end{eqnarray}
for all $n \in \mathbb{Z}$ and $j \in \{1, \ldots, l\}$. The operator $\varpi((\alpha^*)^j)$ restricts on the component $\mathcal{X}_m$ to be the operator 
\begin{eqnarray*}
\varpi((\alpha^*)^j): \bigoplus_{s = 1}^l \mathcal{X}_m^s \rightarrow \bigoplus_{s = 1}^{l-j} \mathcal{X}_{m}^s \oplus \bigoplus_{s = l-j+1}^l \mathcal{X}_{m+1}^s. 
\end{eqnarray*}
which is given in the basis of $l^2(\mathbb{Z}) \otimes \bigoplus_{s = 1}^l l^2(\mathbb{N}_0)$ by 
\begin{eqnarray*}
\varpi((\alpha^*)^j) e_{m+p} \otimes e^s_p &=&  \prod_{r = 1}^j \sqrt{1 - q^{2(lp+s-r)}} e_{m+p} \otimes e^{s+l-j}_{p-1}  \\
&=& \mathcal{S} \prod_{r = 1}^j \sqrt{1 - q^{2(lp+s-r)}}e_{m+p} \otimes e^{s+l-j}_p \5 s \in \{1, \ldots, j\}, \\
\varpi((\alpha^*)^j) e_{m+p} \otimes e^s_p &=& \prod_{r = 1}^j \sqrt{1 - q^{2(lp+s-r)}} e_{m+p} \otimes e^{s-j}_p \5 s \in \{j+1, \ldots, l\}.
\end{eqnarray*}
These operators are independent on the parameter $m \in \mathbb{Z}$. We identify each $\mathcal{X}_m$ with $\bigoplus_{s = 1}^l l^2(\mathbb{N}_0)$ as above, and fix $m \in \mathbb{Z}$ and then view $\varpi((\alpha^*)^j)$ as an operator in $\mathcal{B}(\bigoplus_{s = 1}^l l^2(\mathbb{N}_0))$. Consequently, independently on the choice of $m$, the space $\varpi((\alpha^*)^j \mathcal{L}[n] )$ is now identified with a subspace in $\mathcal{B}(\bigoplus_{s = 1}^l l^2(\mathbb{N}_0))$. If $s \in \{1, \ldots, j\}$, then  the $l^2(\mathbb{N}_0)^s \rightarrow l^2(\mathbb{N}_0)^{s+l-j}$ block of the operator $\varpi((\alpha^*)^j)$ can be written by $\mathcal{S}(1 + T_s)$ where $T_s$ is a compact operator (recall that $q \in (0,1)$), and if $s \in \{j+1, \ldots, l\}$, the  $l^2(\mathbb{N}_0)^s \rightarrow l^2(\mathbb{N}_0)^{s-j}$ block can be written by $(1 + T'_s)$ where $T'_s$ is a compact operator. It follows that $ \varpi( (\alpha^*)^j \mathcal{L}[n] )$ lies in the subspace of $\mathcal{B}(\bigoplus_{s = 1}^l  l^2(\mathbb{N}_0))$ of the operators
\begin{eqnarray}\label{matrix}
\begin{pmatrix} \textbf{0} & X \\ Y & \textbf{0} \end{pmatrix}, \5 X = (\mathcal{K} + \mathbb{C}\mathcal{S}^n)\textbf{1}_{l-j}  \5 Y =  (\mathcal{K} + \mathbb{C}\mathcal{S}^{n+1})\textbf{1}_{j}  
\end{eqnarray}
where $\textbf{1}_j$ is the unit $j \times j$-matrix. Similar analysis shows that the space $\varpi(\alpha^{l-j} \mathcal{L}[n+1])$ is a subspace in \eqref{matrix} as well. This is a complete vector space since it can be identified with a direct sum of copies of $\mathcal{K}$ and $\mathbb{C}$. 

Consider the operator $c\mathcal{S}^{n+1} + A: l^2(\mathbb{N}_0)^s \rightarrow l^2(\mathbb{N}_0)^{s+l-j}$ with $A \in \mathcal{K}$ and $c \in \mathbb{C}$. This is a limit of a Cauchy sequence of operators in $\varpi( (\alpha^*)^j \mathcal{L}[n])$: there is a Cauchy sequence of compact operators $A_i \in \varpi(\mathcal{L}[n])$ which converges to the operator $A$ in the block $l^2(\mathbb{N}_0)^s \rightarrow l^2(\mathbb{N}_0)^s$ and to the zero operator elsewhere, and now
\begin{eqnarray*}
\mathcal{S}(1 + T_s)(1 + \sum_{i = 1}^N (-T_s)^i ) (c \mathcal{S}^n + \mathcal{T} A_N) \in \varpi(\alpha^j) \varpi(\mathcal{L}[n])
\end{eqnarray*}
converges to the operator under the consideration in the limit $N \rightarrow \infty$ because $||T_s|| < 1$. We used the symbol $\mathcal{T}$ for the right inverse of $\mathcal{S}$. Similarly one checks that the operators in the $X$-block of \eqref{matrix} arise as limits of the operators in $\varpi((\alpha^*)^j\mathcal{L}[n])$. So, each element in \eqref{matrix} is a limit of a Cauchy sequence in \eqref{prop3}. It follows that $\overline{\mathbb{C}[SU_{2,q}]^{(nl+j)}}$ is the space \eqref{matrix}. 

Since both modules, $\overline{\mathbb{C}[SU_{2,q}]^{(nl+j)}}$ and \eqref{module2}, are subspaces in $\mathcal{B}(\bigoplus_{s = 1}^l l^2(\mathbb{N}_0))$, we can restrict the linear map
\begin{eqnarray*}
\mathcal{B}(\bigoplus_{s = 1}^l l^2(\mathbb{N}_0)) \times \mathcal{B}(\bigoplus_{s = 1}^l l^2(\mathbb{N}_0)) \rightarrow \mathcal{B}(\bigoplus_{s = 1}^l l^2(\mathbb{N}_0)); \5 (A_1, A_2) \mapsto A_1 A_2^*
\end{eqnarray*}
to define the sesquilinear pairings associated with the $C^*$-module structures. Then $\overline{\mathbb{C}[SU_{2,q}]^{(nl+j)}}$ and \eqref{module2} are isomorphic $C^*$-modules for the algebra $C(\mathbb{WP}_{1,l,q})$ for all $l \in \mathbb{N}$, $n \in \mathbb{Z}$ and $j \in \{1, \ldots, l-1\}$. \5 $\square$ \3

If $n \geq 0$ and $j \in \{1, \ldots, l-1\}$, then the $C(\mathbb{WP}_{1,l,q}) = (\bigoplus_{s = 1}^l \mathcal{K})^+$-module $\overline{\mathbb{C}[SU_{2,q}]^{(nl+j)}}$ is isomorphic to 
\begin{eqnarray*}
((\bigoplus_{s = 1}^l \mathcal{K})^+ \oplus (\bigoplus_{s = 1}^l \mathcal{K})^+) (I_1 \oplus (\bigoplus_{s = 1}^{l-j} P_{n}) \oplus (\bigoplus_{s = l-j+1}^{l} P_{n+1}) ) 
\end{eqnarray*}
and if $n < 0$, then it is isomorphic to 
\begin{eqnarray*}
(\bigoplus_{s = 1}^l \mathcal{K})^+ (1 - (\bigoplus_{s = 1}^{l-j} P_{n}) \oplus (\bigoplus_{s = l-j+1}^{l} P_{n+1})), 
\end{eqnarray*}
see \cite{She14}. So, we get the following theorem. \3

\noindent \textbf{Theorem 1.} The left $C(\mathbb{WP}_{1,l,q}) = (\bigoplus_{s = 1}^l \mathcal{K})^+$-module $\overline{\mathbb{C}[SU_q(2)]^{(nl+j)}}$ is isomorphic to the projective module determined by the projection 
\begin{eqnarray*}
I_1 \oplus (\bigoplus_{s = 1}^{l-j} P_{n}) \oplus (\bigoplus_{s = l-j+1}^{l} P_{n+1})  \5 \text{if} &&\5 n \geq 0, \\
1 - (\bigoplus_{s = 1}^{l-j} P_{n}) \oplus (\bigoplus_{s = l-j+1}^{l} P_{n+1})  \5 \text{if} &&\5 n < 0
\end{eqnarray*}
for all $n \in \mathbb{Z}$ and $j \in \{1, \ldots, l-1\}$.\3

For all $i \in \mathbb{Z}$, the spaces $\overline{\mathbb{C}[SU_{2,q}]^{(i)}}$ are projective modules over the $C^*$-algebra $C(\mathbb{WP}_{1,l,q})$ and all the even K-theory elements can be realized by applying the group structure in K-theory with these modules.\3 

\noindent \textbf{Corollary 1.} The $\mathbb{C}[\mathbb{WP}_{1,l,q}]$-modules $\mathbb{C}[SU_{2,q}]^{(i)}$ and $\mathbb{C}[SU_{2,q}]^{(j)}$ are isomorphic if and only if $i = j$. \3

\noindent Proof. Consider $\mathbb{C}[SU_{2,q}]^{(i)}$ and $\mathbb{C}[SU_{2,q}]^{(j)}$ as subspaces in $\mathcal{B}( \bigoplus_{s = 1}^l l^2(\mathbb{N}_0))$ and as normed
vector spaces with the norms given by the restrictions of the operator norm. Suppose that $U: \mathbb{C}[SU_{2,q}]^{(i)} \rightarrow \mathbb{C}[SU_{2,q}]^{(j)}$ is an isomorphism of $\mathbb{C}[\mathbb{WP}_{1,l,q}]$-modules. Then $U$ induces an isomorphism of normed vector spaces and consequently $U$ extends to an isomorphism of $C(\mathbb{WP}_{1,l,q})$-modules and therefore $i = j$.\5 $\square$ \3

The analysis in 2.3 has been restricted to the case $k = 1$. However, according to \cite{BF12}, the $C^*$-algebra $C(\mathbb{WP}_{k,l,q})$ is isomorphic to $C(\mathbb{WP}_{1,l,q})$ for all coprime positive integers $k$ and $l$. It is expected that the decomposition into the homogeneous components leads to a full analytic realization for the $K$-theory in the cases $k \neq 1$ as well, and that the homogeneous components are mutually non-isomorphic as $\mathbb{C}[\mathbb{WP}_{k,l,q}]$-modules. \3

\noindent \textbf{2.4.} Next we proceed to define a coaction of $B$ on $\mathbb{C}[SU_{2,q}] \otimes M_{\frac{1}{2}}$. We first extend the coaction to the dual Hopf algebra $U_q(\mathfrak{su}_2)$. We will use the same symbol for this coaction 
\begin{eqnarray*}
\theta: U_q(\mathfrak{su}_2) \rightarrow U_q(\mathfrak{su}_2) \otimes B
\end{eqnarray*}
and we require the right $U_q(\mathfrak{su}_2)$-equivariance with respect to the right representation $\partial$: 
\begin{eqnarray}\label{coaction}
\theta(\partial(x) t) = (\partial \otimes \iota)(\theta(x))(\theta t).
\end{eqnarray}
for each $x \in U_q(\mathfrak{su}_2)$ and $t \in \mathbb{C}[SU_{2,q}]$. The coaction will be extended to the tensor products $U_q(\mathfrak{su}_2) \otimes U_q(\mathfrak{su}_2)$ so that
\begin{eqnarray*}
\theta(x \otimes y) = x \otimes y \otimes u^{i+j} 
\end{eqnarray*}
if $x$ is a homogeneous of order $i$ and $y$ is a homogeneous of order $j$. \3

\noindent \textbf{Proposition 4.} There is a unique right $U_q(\mathfrak{su}_2)$-equivariant coaction of $B$ on $U_q(\mathfrak{su}_2)$ which is given on the generators by 
\begin{eqnarray}\label{ad2}
\theta(e) = e \otimes u^{-k-l}, \5 \theta(f) = f \otimes u^{k+l}, \5 \theta(k^{\pm 1}) = k^{\pm 1} \otimes 1
\end{eqnarray}
and extended to the algebra by linearity and 
\begin{eqnarray*}
\theta(x_1 \cdots x_k) = \theta(x_1) \cdots \theta(x_k).
\end{eqnarray*}

\noindent Proof. It is elementary to check that \eqref{ad2} defines a coaction by using the Hopf structure: $\triangle(u) = u \otimes u$ and $\epsilon(u) = 1$. If we study  the elements $\alpha, \alpha^*, \beta, \beta^*$ then a straightforward computation gives that the equivariance \eqref{coaction} holds only if the coaction satisfies \eqref{ad2}. Conversely, \eqref{ad2} together with the product rule determines a  right equivariant coaction on the subspace spanned by these elements. It is therefore sufficient to check that the equivariance remains valid for more general elements $t \in \mathbb{C}[SU_{2,q}]$. 

Notice that
\begin{eqnarray}\label{3}
(\triangle \otimes \iota)(\theta(x)) = \theta(x' \otimes x'')
\end{eqnarray}
for all $x \in U_q(\mathfrak{su}_2)$ since $k$ and $k^{-1}$ in $U_q(\mathfrak{su}_2)$ are coinvariant. 

 We shall proceed by proving that if the equivariance holds on the pair of elements  $t^{\lambda}_{mn}$ and $t^{\lambda'}_{m'n'}$ in $\mathbb{C}[SU_{2,q}]$, then the equivariance holds on the product as well. Take $x \in U_q(\mathfrak{su}_2)$. Then 
\begin{eqnarray*}
\theta(\partial(x) t^{\lambda}_{mn}  t^{\lambda'}_{m'n'}) &=& \theta (\partial(x') t^{\lambda}_{mn}  \partial(x'')t^{\lambda'}_{m'n'}) \\
&=& \theta (\partial(x') t^{\lambda}_{mn} ) \theta (\partial(x'')t^{\lambda'}_{m'n'}) \\
&=& ( \partial \otimes \iota) ( \theta (x')) \theta(t^{\lambda}_{mn} )  ( \partial \otimes \iota)  (\theta (x'')) \theta(t^{\lambda'}_{m'n'}) \\
&=& (\partial \otimes \iota)(\theta(x' \otimes x'')) \theta(t^{\lambda}_{mn}) \theta(t^{\lambda'}_{m'n'}) \\
&=& (\partial \otimes \iota)(\theta(x)) \theta(t^{\lambda}_{mn}t^{\lambda'}_{m'n'}).
\end{eqnarray*}
The first equality holds by Proposition 1, the third holds by the equivariance and the last follows from \eqref{3}. Since $\mathbb{C}[SU_{2,q}]$ is generated by elements which satisfy the equivariance \eqref{coaction}, it is clear that the equivariance holds in the algebra $\mathbb{C}[SU_{2,q}]$. \5 $\square$ \3 
 
We have made $U_q(\mathfrak{su}_2)$ a $B$-comodule algebra. A coaction $\hat{\theta}$ of $B$ on the pre-Hilbert space $\mathbb{C}[SU_{2,q}] \otimes M_{\frac{1}{2}}$ is called a lift of $\theta: \mathbb{C}[SU_{2,q}] \rightarrow \mathbb{C}[SU_{2,q}] \otimes B$, if 
\begin{eqnarray*}
\hat{\theta} =  \theta \otimes \theta' : \mathbb{C}[SU_{2,q}] \otimes M_{\frac{1}{2}} \rightarrow \mathbb{C}[SU_{2,q}] \otimes M_{\frac{1}{2}} \otimes B
\end{eqnarray*}
so that $\theta'$ is a coaction on $M_{\frac{1}{2}}$ (we exploit the notation and assume that the product is applied in the $B$-components arising from the coactions without writing it out explicitly). In addition, we call $\hat{\theta}$ a right $U_q(\mathfrak{su}_2)$-equivariant if 
\begin{eqnarray*}
\hat{\theta}(\partial'(x) (t^{\lambda}_{mn} \otimes v)) = (\partial' \otimes \iota)( \theta(x)) \hat{\theta}(t^{\lambda}_{mn} \otimes v)
\end{eqnarray*} 
for all $x \in U_q(\mathfrak{su}_2)$ and $t^{\lambda}_{mn} \otimes v \in \mathbb{C}[SU_{2,q}] \otimes M_{\frac{1}{2}}$. \3

\noindent \textbf{Proposition 5.} The right $U_q(\mathfrak{su}_2)$-equivariant lifted coactions of $\theta$ on the  pre-Hilbert space $\mathbb{C}[SU_{2,q}] \otimes M_{\frac{1}{2}}$ are parametrized by $i \in \mathbb{Z}$, and for all $i \in \mathbb{Z}$, the lifted coaction $\hat{\theta}_i$ is determined by   
\begin{eqnarray*}
\hat{\theta}_i(t^{\lambda}_{mn} \otimes v) =  \theta(t^{\lambda}_{mn}) \otimes \theta_i(v)  
\end{eqnarray*}
where $\theta_i: M_{\frac{1}{2}} \rightarrow M_{\frac{1}{2}} \otimes B$ is the coaction of $B$ determined by  
\begin{eqnarray*}
\theta_i(e_+) = e_+ \otimes u^i, \5 \theta_i(e_-) = e_- \otimes u^{i + k + l}.
\end{eqnarray*}
The lifted coactions $\hat{\theta}_i$ are mutually nonequivalent. \3

\noindent Proof. It is straightforward to verify that $\hat{\theta}_i$ and $\theta_i$ are coactions of $B$ for each $i \in \mathbb{Z}$ and that $\theta_i$ are all the possible coactions on $M_{\frac{1}{2}}$ which are $U_q(\mathfrak{su}_2)$-equivariant under the representation $\varrho_{\frac{1}{2}}$: 
\begin{eqnarray*}
\theta_i(\varrho_{\frac{1}{2}}(x) v) = (\varrho_{\frac{1}{2}}(x) \otimes \iota)(\theta_i(x))(\theta(v)).
\end{eqnarray*}
Then we use the equivariance of the coations $\theta$ and $\theta_i$, and \eqref{3} to obtain 
\begin{eqnarray*}
\hat{\theta}_i(\partial'(x) t^{\lambda}_{mn} \otimes v) &=& \hat{\theta}_i(\partial(x') t^{\lambda}_{mn} \otimes \varrho_{\frac{1}{2}}(x'') v) \\
&=& \theta(\partial(x') t^{\lambda}_{mn}) \otimes \theta_i( \varrho_{\frac{1}{2}}(x'') v) \\
&=& (\partial \otimes \iota )( \theta(x') )\theta(t^{\lambda}_{mn}) \otimes (\varrho_{\frac{1}{2}} \otimes \iota)(\theta(x''))\theta_i(v) \\
&=& (\partial \otimes \varrho_{\frac{1}{2}} \otimes \iota) (\theta(x') \otimes \theta(x''))  (\theta(t^{\lambda}_{mn}) \otimes \theta_i(v)) \\
&=& (\partial \otimes \varrho_{\frac{1}{2}} \otimes \iota) (\triangle \otimes \iota) (\theta(x)) (\theta(t^{\lambda}_{mn}) \otimes \theta_i(v)) \\
&=& (\partial' \otimes \iota) (\theta(x)) (\hat{\theta}_i(t^{\lambda}_{mn} \otimes v)).
\end{eqnarray*}
So, the  coaction $\hat{\theta}_i$ is right equivariant. 

Conversely, since $1$ is the cyclic vector for the GNS representation, a lifted coaction is fully determined by its coaction on to the subspace $1 \otimes M_{\frac{1}{2}}$ which is an irreducible representation of $U_q(\mathfrak{su}_2)$ of the highest weight $\frac{1}{2}$. The right equivariance implies that $\hat{\theta}$ restricts on this subspace to $\iota \otimes \theta_i$ for some $i \in \mathbb{Z}$. So $\hat{\theta}_i$ are all the possible solutions. An equivalence of coactions of $B$ sends homogeneous elements to each other and preserves the degree of homogeneity. So, if $i \neq j$ in $\mathbb{Z}$, then $\theta_i$ and $\theta_j$ are nonequivalent coactions. It follows that $\hat{\theta}_i$ and $\hat{\theta}_j$ are nonequivalent. \5 $\square$ \3

\noindent \textbf{2.5.} We shall build a spectral triple on a coinvariant subspace of the Hilbert space completion of $\mathbb{C}[SU_{2,q}] \otimes M_{\frac{1}{2}}$. The goal is to construct an explicit model and for this reason we fix the lifted coaction: 
\begin{eqnarray*}
\hat{\theta} := \hat{\theta}_{-k} = \theta \otimes \theta_{-k}
\end{eqnarray*}
where $k$ is an arbitrary fixed integer. Now the weight vectors $e_+$ and $e_-$ in $M_{\frac{1}{2}}$ are homogeneous of orders $-k$ and $l$ under $\theta_{-k}$. \3

\noindent \textbf{Proposition 6.} The coinvariant subspaces in the pre-Hilbert spaces $\mathbb{C}[SU_{2,q}]$ and $\mathbb{C}[SU_{2,q}] \otimes M_{\frac{1}{2}}$  are linearly spanned by the sets of vectors
\begin{eqnarray*}
&& t^{\lambda}_{p(l+k),p(l-k)} \5 \text{and} \\
&&  t^{\lambda}_{p(l+k)-\frac{1}{2},p(l-k) - \frac{1}{2}} \otimes e_+ , \5  t^{\lambda}_{p(l+k) - \frac{1}{2},p(l-k) + \frac{1}{2}} \otimes e_- 
\end{eqnarray*}
respectively. The index $\lambda$ runs over all the highest weights, $\lambda \in \frac{1}{2} \mathbb{N}_0$, and $p$ gets all the half-integer values for which these vectors are defined. \3

\noindent Proof. Recall that $m$ and $n$ are both integers, or both half-integers. In the notation of Proposition 2, we need to find all the pairs $(m,n)$ which solve $(m+n)k = (m-n)l$. Since $k$ and $l$ are coprime, the solutions have to satisfy $m+n = xl$ and $m-n = xk$ for some $x \in \mathbb{Q}$. On the other hand, $m+n$ and $m-n$ are always integers, and since $l$ and $k$ are coprime, $x$ has to be an integer. It follows that the solutions are $m = p(l+k)$ and $n = p(l-k)$ with $p \in \frac{1}{2} \mathbb{Z}$. Thus the coinvariant subspace in $\mathbb{C}[SU_{2,q}]$  linearly spanned by 
\begin{eqnarray*}
 t^{\lambda}_{p(l+k),p(l-k)}
\end{eqnarray*}
for all $\lambda \in \frac{1}{2} \mathbb{N}_0$ and $p$ such that the vectors are well defined. 

Since $e_+$ and $e_-$ are homogeneous elements of orders $-k$ and $l$, it is sufficient to find the homogeneous subspaces of orders $k$ and $-l$ in the algebra $\mathbb{C}[SU_{2,q}]$ to get the coinvariant subspace in $\mathbb{C}[SU_{2,q}] \otimes M_{\frac{1}{2}}$. The same analysis as above gives the subspaces spanned by the vectors 
\begin{eqnarray*}
 t^{\lambda}_{p(l+k)-\frac{1}{2},p(l-k)-\frac{1}{2}} \5 \text{and} \5 t^{\lambda}_{p(l+k) - \frac{1}{2},p(l-k) + \frac{1}{2}}
\end{eqnarray*}
for all $\lambda \in \frac{1}{2} \mathbb{N}_0$ and $p$ such that the vectors are well defined. \5 $\square$ \3

\noindent \textbf{Lemma 2.} The linear basis of the coinvariant subspace of $\mathbb{C}[SU_{2,q}] \otimes M_{\frac{1}{2}}$ is given in the notation $|j m \mu \updownarrows \ra$ by
\begin{eqnarray*}
  |j, p(l+k) - \frac{1}{2}, p(l-k) \uparrow \ra \5 \text{and} \5 |j, p(l+k) - \frac{1}{2}, p(l-k) \downarrow \ra. 
\end{eqnarray*}
The index $j$ runs over $\frac{1}{2}\mathbb{N}_0$, and for any fixed $j$ the index $p \in \frac{1}{2} \mathbb{Z}$ runs over the values for which $p(l+k) \in \{-j, \ldots,j, j+1\}$ in the subspace labeled by $\uparrow$, and over the values for which $p(l+k) \in \{-j+1, \ldots, j-1,j\}$ in the subspace labeled by $\downarrow$. \3

\noindent Proof. Using the Clebsch-Gordan coefficients introduced in 1.4 we find that
\begin{eqnarray*}
& &q^{-p(l+k)+\frac{1}{2}}([2j^- + 1]_q)^{-\frac{1}{2}} |j-1, p(l+k) - \frac{1}{2}, p(l-k) \uparrow \ra \\ 
&=& -S_{j, p(l-k)} t^{j^-}_{p(l+k) - \frac{1}{2}, p(l-k) + \frac{1}{2}} \otimes e_- + C_{j,p(l-k)} t^{j^-}_{p(l+k) - \frac{1}{2}, p(l-k) - \frac{1}{2}} \otimes e_+
\end{eqnarray*}
if $j \geq 1$ and $p \in \frac{1}{2} \mathbb{Z}$ is so that $p(l+k) - \frac{1}{2} \in \{-(j-1)^+, \ldots,(j-1)^+-1, (j-1)^+\}$ and $p(l-k) \in \{-(j-1), \ldots, j-2, j-1\}$ (recall the parametrization of the basis $|j m \mu \updownarrows \ra$ in 1.4). The first condition is equivalent to $p(l+k) \in \{-j+1, \ldots, j-1, j\}$ and the second one follows from this condition. Similarly,  
\begin{eqnarray*}
& &q^{-p(l+k)+\frac{1}{2}}([2j^- + 1]_q)^{-\frac{1}{2}} |j, p(l+k) - \frac{1}{2}, p(l - k) \downarrow \ra \\
&=& C_{j,p(l-k)} t^{j^-}_{p(l+k) - \frac{1}{2}, p(l-k) + \frac{1}{2}} \otimes e_- + S_{j,p(l-k)} t^{j^-}_{p(l+k) - \frac{1}{2}, p(l-k) - \frac{1}{2}} \otimes e_+
\end{eqnarray*}
if $j \geq \frac{1}{2}$ and $p \in \frac{1}{2} \mathbb{Z}$ is so that $p(l+k) \in \{-j+1, \ldots, j-1,j\}$ holds. Therefore, the vectors listed in Lemma 2 form a basis for the pre-Hilbert space of coinvariant vectors. \5 $\square$ \3 

For what follows we need to compute dimensions of certain subspaces in the coinvariant pre-Hilbert spaces of Proposition 6. For all $j \in \frac{1}{2} \mathbb{N}_0$, denote by $V_j^{\updownarrows}$ the subspace of $\mathbb{C}[SU_{2,q}] \otimes M_{\frac{1}{2}}$ spanned by the vectors of Lemma 2 with the fixed indexes $j$ and $\updownarrows$, and $p \in \frac{1}{2} \mathbb{Z}$ runs over all its possible values. \3

\noindent \textbf{Proposition 7.} For each $j \in \frac{1}{2} \mathbb{N}_0$, the relation $\text{dim}(V_{j}^{\uparrow}) = \text{dim}(V_{j+1}^{\downarrow})$ holds, $\text{dim}(V_0^{\downarrow}) = 0$ and 
\begin{quote}
\textbf{i.} if $k+l$ is even, then $V_j^{\updownarrows}$ are nonzero only if $j$ is an integer, and if $j \in \mathbb{N}$,  
\begin{eqnarray*}
\text{dim}(V_j^{\downarrow}) =  \Big[ \frac{j}{\frac{1}{2}(l+k)} \Big] + \Big[ \frac{j - 1}{\frac{1}{2}(l+k)} \Big] + 1,
\end{eqnarray*}
\textbf{ii.} if $l+k$ is odd, then 
\begin{eqnarray*}
\text{dim}(V_j^{\downarrow}) =  \Big[ \frac{j}{l+k} \Big] + \Big[ \frac{j-1}{l+k} \Big]  + 1&  \5 &\text{if $j \in \mathbb{N}$} \\
\text{dim}(V_j^{\downarrow}) = \Big[ \frac{j}{l+k} + \frac{1}{2} \Big] + \Big[ \frac{j-1}{l+k} + \frac{1}{2} \Big]&  \5 &\text{if $j \in \mathbb{N}_0 + \frac{1}{2}$}.
\end{eqnarray*}
$[ a ]$ denotes the integer part of $a$. 
\end{quote}

\noindent Proof. It follows from Lemma 2 that for any $j \in \frac{1}{2} \mathbb{N}_0$ we have $\text{dim}(V_{j}^{\uparrow}) = \text{dim}(V_{j+1}^{\downarrow})$. Thus, it is sufficient to find the dimensions of $V_j^{\downarrow}$. By definition, $\text{dim}(V_0^{\downarrow}) = 0$.

Let $k+l$ be even. Now $p(l+k)$ is an integer for all $p \in \frac{1}{2} \mathbb{Z}$ and so $j$ must be an integer. We need to find all $p \in \frac{1}{2}\mathbb{Z}$ so that $p(l+k) \in \{-j+1,-j+2, \ldots, j\}$. For all $j \in \mathbb{N}$, there is the solution $p = 0$ and if $j < \frac{1}{2}(l+k)$, then this is the only solution. If $j = \frac{1}{2}(l+k)$, then $p$ can get two values, $0$ and $\frac{1}{2}$. If $ \frac{1}{2}(l+k) + 1 \leq j < l+k$, then $p$ can get three values, $0$ and $\pm \frac{1}{2}$. At $j = l+k$ we get the fourth solution $p = 1$. Continuing like this we observe that a  new pair of solutions for $p$ appear  whenever $j$ grows by $\frac{1}{2}(l+k)$ and therefore i follows. 

Let $k+l$ be odd. Now there are two cases: if $j$ is an integer, then the weights $p(k+l)$ are integers and so $p$ has to be an integer. This case  is essentially the same as above except that  a new pair of weights appear whenever $j$ is shifted by $(l+k)$. This gives the first part of ii. If $j \in \mathbb{N}_0 + \frac{1}{2}$, then $p$ has solutions for $p \in \mathbb{Z} + \frac{1}{2}$. The first $j \in \mathbb{N}_0 + \frac{1}{2}$ for which $V^{\downarrow}_j$ is nonzero is $j = \frac{1}{2}(l+k)$: in this case $p =  \frac{1}{2}$. Then at $j =  \frac{1}{2}(l+k) + 1$ we get the second solution $p = -\frac{1}{2}$. Now, whenever $j$ is shifted by $(l+k)$ one gets another pair of solutions for $p$ and so 
\begin{eqnarray*}
\text{dim}(V_j^{\downarrow}) = \Big[ \frac{j + \frac{l+k}{2}}{l+k} \Big] + \Big[ \frac{j + \frac{l+k}{2} - 1}{l+k} \Big] =  \Big[ \frac{j}{l+k} + \frac{1}{2} \Big] + \Big[ \frac{j-1}{l+k} + \frac{1}{2} \Big].
\end{eqnarray*}
The second part of ii follows. \5 $\square$ \3

For all $\lambda \in \frac{1}{2} \mathbb{N}_0$, denote by $V_{\lambda}$ the subspace of $\mathbb{C}[SU_{2,q}]$ which is spanned by the vectors $ t^{\lambda}_{p(l+k),p(l-k)}$ for all $p \in \frac{1}{2} \mathbb{Z}$ for which $p(l+k) \in \{-\lambda, \ldots, \lambda\}$. Similar analysis of the weight parameters can be applied to prove the following. \3 

\noindent \textbf{Proposition 8.} 
\begin{quote}
\textbf{i.} If $k+l$ is even, then $V_{\lambda}$ is nonzero only if $\lambda$ is an integer, and for all $\lambda \in \mathbb{N}_0$:  
\begin{eqnarray*}
\text{dim}(V_{\lambda}) = 2 \Big[ \frac{\lambda}{\frac{1}{2}(l+k)} \Big] + 1.
\end{eqnarray*}
\textbf{ii.} If $l+k$ is odd then 
\begin{eqnarray*}
\text{dim}(V_{\lambda})   =  2 \Big[ \frac{\lambda}{l+k} \Big]  + 1&  \5 &\text{in the case $\lambda$ is an integer} \\
\text{dim}(V_{\lambda}) = 2 \Big[ \frac{\lambda}{l+k} + \frac{1}{2} \Big]&  \5 &\text{otherwise.}
\end{eqnarray*}
\end{quote}

\section{Dirac Operator Quantum Weighted Projective Spaces}

In this section we construct two models for Dirac operators on coinvariant Hilbert spaces. The first one in 3.1 is based on the restriction of the Dirac operator on $SU_{2,q}$ on the coinvariant component of $\hil$. There is no natural chiral grading in the $SU_{2,q}$ model and this construction leads to an odd spectral triple.  An even spectral triple (i.e. chirally graded) is developed in 3.2. The Fredholm modules arising from these spectral triples represent the trivial class in K-homology. \3 

\noindent \textbf{3.1.} We use the coaction $\hat{\theta} = \theta \otimes \theta_{-k}$ on $\mathbb{C}[SU_{2,q}] \otimes M_{\frac{1}{2}}$. In 2.5 we have a complete description of the pre-Hilbert space of coinvariant spinors in $\mathbb{C}[SU_{2,q}] \otimes M_{\frac{1}{2}}$ and we denote $\hil^{\theta}$ its Hilbert space completion with respect to the Haar functional. The restriction of the representation $\pi:= \pi_h \otimes \iota$ on the coinvariant subalgebra $\mathbb{C}[\mathbb{WP}_{k,l,q}]$ determines a $*$-representation $\pi: \mathbb{C}[\mathbb{WP}_{k,l,q}] \rightarrow \mathcal{B}(\hil^{\theta})$. 

We define the Dirac operator on $\hil^{\theta}$ by restricting the Dirac operator $D$ introduced in 1.5 to $\hil^{\theta}$. In addition, it is convenient to shift the spectrum so that the positive and negative eigenvalues appear pairwise and their multiplicities match. So, let $D^{\theta}$ denote the densely defined self-adjoint operator
\begin{eqnarray*}
D^{\theta} = D + \frac{1}{2}
\end{eqnarray*}
on $\hil^{\theta}$. The restriction is well defined since by Lemma 2, $\hil^{\theta}$ has a basis consisting of the eigenvectors of the Dirac operator on $SU_{2,q}$: $D^{\theta}$ has the eigenvalue $2(j+1)$ on $V^{\uparrow}_j$, and the eigenvalue $-2(j+1)$ on $V^{\downarrow}_{j+1}$. The dimensions are listed in Proposition 7 and they depend on $l$ and $k$. In particular, we always have $V^{\downarrow}_0 = V^{\downarrow}_{\frac{1}{2}} = 0$ and $\text{dim}(V^{\uparrow}_j) = \text{dim}(V^{\downarrow}_{j+1})$. \3

\noindent \noindent \textbf{Theorem 2.} The data $(\mathbb{C}[\mathbb{WP}_{k,l,q}], \pi, \hil^{\theta}, D^{\theta})$ define a $2$-summable odd spectral triple over the quantum projective plane $\mathbb{C}[\mathbb{WP}_{k,l,q}]$ for all positive coprime integers $k$ and $l$. \3

\noindent Proof. Since this model is based on a restriction of a spectral triple on $SU_{2,q}$ most of the details are automatic. All the commutators $[D^{\theta}, \pi(t)]$ extend to bounded operators for all $t \in \mathbb{C}[\mathbb{WP}_{k,l,q}]$ since this is the case for the quantum group. Also this representation is faithful for the same reason. For the summability, we shall consider the case where $k+l$ is even. The sum of eigenvalues of $(D^{\theta})^{-2}$ over the subspaces $V^{\updownarrows}_{j}$ for $0 \leq j \leq N$ gives 
\begin{eqnarray*}
\sigma_N((D^{\theta})^{-2}) = 2\sum_{0 \leq j \leq N} ( \Big[ \frac{j}{\frac{1}{2}(k+l)} \Big] + \Big[ \frac{j-1}{\frac{1}{2}(k+l)} \Big] + 1) (2(j+1))^{-2}.
\end{eqnarray*}
In the limit $N \rightarrow \infty$, $\sigma_N$ diverges logarithmically. Thus the summability is $2$. The same argument gives the same summability if $l+k$ is odd. \5 $\square$ \3

Associated to the odd spectral triple of Theorem 2 there is an odd Fredholm module over the $C^*$-algebra $C(\mathbb{WP}_{k,l,q})$: the GNS representation defines the action of $C(\mathbb{WP}_{k,l,q})$ on the Hilbert space $\hil^{\theta}$, and the sign operator $F^{\theta} = D^{\theta}|D^{\theta}|^{-1}$ is the bounded Fredholm operator of the Fredholm module. Since the algebra $C(\mathbb{WP}_{k,l,q})$ is isomorphic to $(\bigoplus_{s = 1}^l \mathcal{K})^+$ the odd K-homology of $C(\mathbb{WP}_{k,l,q})$ is trivial. So, the K-homology class of this Fredholm module represents the trivial element. \3

\noindent \textbf{3.2.} Denote by $\hil'$ the Hilbert space completion of the coinvariant component $\mathbb{C}[SU_{2,q}]^{(0)}$ with respect to the Haar measure. The Hilbert space $\hil'$ is a representation space for $\mathbb{C}[\mathbb{WP}_{k,l,q}]$ under $\pi_h$. We define the Hilbert space of spinors to be the Hilbert space direct sum $\hil' \oplus \hil'$ and let $\mathbb{C}[\mathbb{WP}_{k,l,q}]$ act on it under $\pi_h \oplus \pi_h$. Let us fix the basis consisting of
\begin{eqnarray*}
t^{\lambda, \uparrow}_{p(l+k), p(l-k)} \5 \text{and} \5 t^{\lambda, \downarrow}_{p(l+k), p(l-k)}
\end{eqnarray*}
where we have used the upper-scripts $\uparrow$ and $\downarrow$ to label the first and the second component of the direct sum. Then we define the densely defined self-adjoint operator by
\begin{eqnarray*}
D' t^{\lambda, \uparrow}_{p(l+k), p(l-k)} = (\lambda  + 1) t^{\lambda, \downarrow}_{p(l+k), p(l-k)},\5  D' t^{\lambda, \downarrow}_{p(l+k), p(l-k)} = -(\lambda  + 1) t^{\lambda, \uparrow}_{p(l+k), p(l-k)}.  
\end{eqnarray*}
The spectrum is fixed so that for $k = l = 1$, when the algebra $\mathbb{C}[\mathbb{WP}_{k,l,q}]$ is a deformation of the coordinate algebra of a sphere, the spectrum coincides with the classical Dirac spectrum \cite{GVF01}. The chirality operator is defined to be the diagonal operator 
\begin{eqnarray*}
\omega t^{\lambda, \uparrow}_{p(l+k), p(l-k)} = t^{\lambda, \uparrow}_{p(l+k), p(l-k)} \5  \omega t^{\lambda, \downarrow}_{p(l+k), p(l-k)} = - t^{\lambda, \downarrow}_{p(l+k), p(l-k)}.
\end{eqnarray*}
 
\noindent \textbf{Theorem 3.} The data $(\mathbb{C}[\mathbb{WP}_{k,l,q}], \pi_h \oplus \pi_h, \hil' \oplus \hil', D', \omega)$ defines a $2$-summable even spectral triple over the quantum projective plane for all positive coprime integers $k$ and $l$.\3

\noindent Proof. We need to check that the commutators $[D', (\pi_h \oplus \pi_h)(t)]$ extend to bounded operators for all $t \in \mathbb{C}[\mathbb{WP}_{k,l,q}]$.  First we consider the Hilbert space completion of $\mathbb{C}[SU_{2,q}]$ which is denoted by $L^2(SU_{2,q})$. Then we take the Hilbert space sum $L^2(SU_{2,q}) \oplus L^2(SU_{2,q})$. We label by $\uparrow$ and $\downarrow$ the basis vectors in the first and second summand, and define a densely defined auxiliary Dirac operator on $L^2(SU_{2,q}) \oplus L^2(SU_{2,q})$ by 
\begin{eqnarray*}
Q \eta(t^{\lambda, \uparrow}_{mn}) = (\lambda  + 1) \eta(t^{\lambda, \downarrow}_{mn}) \5 \text{and} \5 Q \eta(t^{\lambda, \downarrow}_{mn}) = -(\lambda  + 1)\eta(t^{\lambda, \uparrow}_{mn}).
\end{eqnarray*}
The algebra $\mathbb{C}[SU_{q,2}]$ acts on $L^2(SU_{2,q}) \oplus L^2(SU_{2,q})$ under the direct sum of GNS representations $\pi_h \oplus \pi_h$. Next we check that the commutators $[Q, (\pi_h \oplus \pi_h)(t)]$ extend to bounded operators for all $t \in \mathbb{C}[SU_{2,q}]$. It is clearly sufficient to check this for the generators. For $t = \alpha$ we get the following action on the basis
\begin{eqnarray*}
\pi_h(\alpha) \eta(t^{\lambda, \updownarrows}_{mn}) &=& \eta(t^{\frac{1}{2}}_{\frac{1}{2} \frac{1}{2}} t^{\lambda, \updownarrows}_{mn}) \\
&=& \sum_{\mu = \lambda-\frac{1}{2}}^{\lambda +\frac{1}{2}}C_q  \begin{pmatrix} \frac{1}{2} & \lambda & \mu \\ \frac{1}{2}  & m & m + \frac{1}{2}  \end{pmatrix} C_q  \begin{pmatrix} \frac{1}{2}  & \lambda & \mu \\ \frac{1}{2}  & n & n+\frac{1}{2}  \end{pmatrix} \eta(t^{\mu, \updownarrows}_{m+\frac{1}{2} ,n+\frac{1}{2} }).
\end{eqnarray*}
Then a straightforward computation gives 
\begin{eqnarray*}
[Q, \pi_h(\alpha)] \eta(t^{\lambda, \updownarrows}_{mn}) &=& - \frac{1}{2}C_q  \begin{pmatrix} \frac{1}{2} & \lambda & \lambda - \frac{1}{2} \\ \frac{1}{2}  & m & m + \frac{1}{2}  \end{pmatrix} C_q  \begin{pmatrix} \frac{1}{2}  & \lambda & \lambda - \frac{1}{2} \\ \frac{1}{2}  & n & n+\frac{1}{2}  \end{pmatrix}   \eta(t^{\lambda - \frac{1}{2}, \downuparrows}_{m+\frac{1}{2} ,n+\frac{1}{2} }) \\
 & & +  \frac{1}{2}C_q  \begin{pmatrix} \frac{1}{2} & \lambda & \lambda + \frac{1}{2} \\ \frac{1}{2}  & m & m + \frac{1}{2}  \end{pmatrix} C_q  \begin{pmatrix} \frac{1}{2}  & \lambda & \lambda + \frac{1}{2} \\ \frac{1}{2}  & n & n+\frac{1}{2}  \end{pmatrix}   \eta(t^{\lambda + \frac{1}{2}, \downuparrows}_{m+\frac{1}{2} ,n+\frac{1}{2} }) 
\end{eqnarray*}
The Clebsch-Gordan coefficients are bounded in $\lambda, m$ and $n$ and therefore $[Q, (\pi_h \oplus \pi_h)(\alpha)]$ extends to a bounded operator. Similar computation proves that the same holds for $[Q, (\pi_h \oplus \pi_h)(\beta)]$. It follows that $[Q, (\pi_h \oplus \pi_h)(t)]$ extends to a bounded operator for all $t \in \mathbb{C}[SU_{2,q}]$. Now the coinvariant spectral triple $(\mathbb{C}[\mathbb{WP}_{k,l,q}], \pi_h \oplus \pi_h, \hil' \oplus \hil', D', \omega)$ can be reduced from the auxiliary model by restricting $L^2(SU_{2,q})$ to the coinvariant component  and by restricting $\mathbb{C}[SU_{2,q}]$ to $\mathbb{C}[\mathbb{WP}_{k,l,q}]$ and $Q$ to $D'$. Then $[D', \pi_h \oplus \pi_h(t)]$ extend to a bounded operator on $\hil' \oplus \hil'$ for all $t \in \mathbb{C}[\mathbb{WP}_{k,l,q}]$. 

The summability is computed exactly as in Theorem 2 with the help of Proposition 8. The representation $\pi_h \oplus \pi_h$ is faithful since both components are faithful. This is an even spectral triple since there is the chirality operator $\omega$ satisfying $\{\omega, D' \} = 0$, $\omega^2 = 1$, $\omega^* = \omega$ and $[(\pi_h \oplus \pi_h)(t), \omega] = 0$ for all $t \in \mathbb{C}[\mathbb{WP}_{k,l,q}]$.\5 $\square$\3

The even K-homology groups of $C(\mathbb{WP}_{k,l,q})$ are isomorphic to $\mathbb{Z}^{\oplus l + 1}$. There is an even Fredholm module arising from the spectral triple of Theorem 3: the $C^*$-algebra $C(\mathbb{WP}_{k,l,q})$ acts on $\hil' \otimes \hil'$ under the sum of two GNS representation, the chirality operator  $\omega$ is defined above, and the bounded Fredholm operator of this Fredholm module is the operator
\begin{eqnarray*}
F' t^{\lambda, \uparrow}_{p(l+k), p(l-k)} =t^{\lambda, \downarrow}_{p(l+k), p(l-k)} \5 \text{and} \5 F' t^{\lambda, \downarrow}_{p(l+k), p(l-k)} = - t^{\lambda, \uparrow}_{p(l+k), p(l-k)}
\end{eqnarray*}
This Fredholm module is degenerate: the operator $F'$ commutes with the representation $\pi_h \oplus \pi_h$, $(F')^* = F'$ and $(F')^2 = 1$. In particular, the K-homology class is zero.

We can construct more even spectral triples by employing the homogeneous components of \eqref{module} in the definition of the Hilbert space. For example, we can take $\hil'_n$ to be the completion of $\mathbb{C}[SU_{2,q}]^{(n)}$ and then apply the Hilbert space $\hil'_n \oplus \hil'_n$. This is a representation space for $\mathbb{C}[\mathbb{WP}_{k,l,q}]$ under $\pi_h \oplus \pi_h$, and the chirality operator can be defined as above. The Dirac operator of Theorem 3 can be specialized to this case in the obvious way; the proof for the boundedness of the commutators of the Dirac operator and the representation operators goes exactly as in Theorem 3. Moreover, according to 2.3, the spectral triples obtained this way are mutually nonequivalent for all $n \in \mathbb{Z}$, at least in the case of a quantum teardrops. The associated Fredholm modules are still degenerate for all $n \in \mathbb{Z}$. One may also find spectral triples based on the Hilbert space $\hil'_n \oplus \hil'_m$ for some $m \neq n$. Such models might determine nontrivial K-homology classes. \3 

\noindent \textbf{3.3.} Let us study the case $l = k = 1$ with more details. The algebra $\mathbb{C}[\mathbb{WP}_{1,1,q}]$ is isomorphic to a 2-dimensional quantum sphere algebra. Now $(l-k) = 0$. The Hilbert space $\hil^{\theta}$ is spanned by the vectors $|j m 0 \uparrow \rangle$ and $|j m 0 \downarrow \rangle$ where $j \in \mathbb{N}_0$ and $m$ gets all the possible values for given $j$. So, the index $m$ runs over all the weight vectors for the left $U_q(\mathfrak{su}_2)$ representation. Consequently, there is a remaining $U_q(\mathfrak{su}_2)$ symmetry from the left. The right symmetry is missing. The dimensions of the Dirac eigenspaces are given by $\text{dim}(V_j^{\uparrow}) = \text{dim}(V_{j+1}^{\downarrow}) = 2j + 1$.  

Now consider the Hilbert space $\hil' \oplus \hil'$. The basis vectors $t^{\lambda, \updownarrows}_{p 0}$ are labeled by $\lambda \in \mathbb{N}_0$ and $p \in \mathbb{Z}$ so that $-\lambda \leq p \leq \lambda$. Again, the left $U_q(\mathfrak{su}_2)$ symmetry is there but the right is missing. The Dirac spectrum now coincides with the usual Dirac spectrum on a $2$-sphere: 
\begin{eqnarray*}
\text{spec}(D') = \{ \pm (\lambda + 1): \lambda \in \mathbb{N}_0 \}
\end{eqnarray*}
and the eigenvalues $\pm (\lambda + 1)$ have the multiplicity $2\lambda + 1$. 

The case $k = 1$ and $l > 1$ corresponds to the standard tear drop with a $\mathbb{Z}_l$-isotropy. Now the left and the right $U_q(\mathfrak{su}_2)$ symmetries are lost since no longer all weight spaces occur.

\bibliographystyle{plain}

\begin{thebibliography}{00}
\bibitem{BF12} Brzezinski T., Fairfax S. A.: Quantum Teardrops, Commun. Math. Phys. 316, 151-170 (2012)
\bibitem{DLSSV05} Dabrowski L., Landi G., Sitarz A., van Suijlekom W., Varilly J.C.: The Dirac Operator on $SU_q(2)$, Commun. Math. Phys. 259, 729-759 (2005)
\bibitem{Dri86} Drinfeld V. G.: Quantum groups; in: Gleason A. M.: Proc. I.C.M. Berkeley (1986) 
\bibitem{GVF01} Gracia-Bondia J. M., Varilly J. C., Figueroa H.: Elements of Noncommutative Geometry (Birkhauser 2001)
\bibitem{Har13a} Harju A., J.: Spectral Triples on Proper Etale Groupoids, to appear in J. Noncom. Geom, preprint: arXiv (2014) 
\bibitem{Jim85} Jimbo M.: A $q$-Difference Analogue of $U(g)$ and the Yang-Baxter Equation, Lett. Math. Phys. 10, 63-69 (1985)
\bibitem{KS97} Klimyk A. ,Schmudgen K.: Quantum groups and their representations (Springer-Verlag 1998)
\bibitem{Moe02} Moerdijk I.: Orbifolds as groupoids: an introduction. Orbifolds in mathematics and physics (Madison, WI, 2001), Contemp. Math. 310, 205-222 (Amer. Math. Soc., Providence, RI, 2002)
\bibitem{MM03} Moerdijk I.,  Mrcun J.: Introduction to Foliations and Lie Groupoids (Cambridge University Press 2003)
\bibitem{NT10} Neshveyev S., Tuset L.: The Dirac operator on compact quantum groups, J. Reine Angew. Math. 641, 1-20 (2010)
\bibitem{NT11}  Neshveyev S., Tuset L.: Notes on the Kazhdan–Lusztig Theorem on Equivalence of the Drinfeld Category and the Category of $U_q \mathfrak{g}$-Modules, Algebr. Represent. Theory 14, 897-948 (2011)
\bibitem{RV08} Rennie, A. C., Varilly, J. C.: Orbifolds are not commutative geometries, J. Australian Math. Soc. 84, 109-116 (2008)
\bibitem{Sat56} Satake I.: On a generalization of the notion of manifold,  Proc. Nat. Acad. Sci. U.S.A. 42, 359-363 (1956)
\bibitem{SV13} Sitarz A., Venselaar J. J.: Real spectral triples on 3-dimensional noncommutative lens spaces, preprint: arXiv (2013)
\bibitem{She14} Sheu A. J. L.: The Structure of Line Bundles over Quantum Teardrops, SIGMA 10 (2014)
\bibitem{Sle85} Slebarski S.: Dirac operators on a compact Lie group, Bull. London Math. Soc. 17, 579-583 (1985) 
\bibitem{Wor87b} Woronowicz S. L.: Twisted SU(2) Group. An Example of a Non-Commutative Differential Calculus, Publ. Res. Inst. Math. Sci. 23, 117-181 (1987)




\end{thebibliography}

\end{document}